# MULTITYPE RANDOMIZED REED–FROST EPIDEMICS AND EPIDEMICS UPON RANDOM GRAPHS

By Peter Neal

*University of Manchester*

We consider a multitype epidemic model which is a natural extension of the randomized Reed–Frost epidemic model. The main result is the derivation of an asymptotic Gaussian limit theorem for the final size of the epidemic. The method of proof is simpler, and more direct, than is used for similar results elsewhere in the epidemics literature. In particular, the results are specialized to epidemics upon extensions of the Bernoulli random graph.

**1. Introduction.** The randomized Reed–Frost epidemic is a very general model for homogeneously mixing SIR (susceptible → infective → removed) epidemic models; see [14] and [18]. That is, while infectious, an infective, $i$ say, has probability $V_i$ of making an infectious contact with any given susceptible member of the population. The random variables $\{V_i\}$ are assumed to be independent and identically distributed. Furthermore, the only transitions in state are: from susceptible to infective, and from infective to removed. Since infectious contacts with nonsusceptible individuals have no effect, we do not need to make assumptions about the relationship between an infective and a nonsusceptible individual.

It is trivial to show that the generalized stochastic epidemic (see [5]), where individuals have independent and identically distributed infectious lifetimes, and while infectious make (potential) infectious contacts at the points of a homogeneous Poisson point process, is a special case of the randomized Reed–Frost epidemic. It was shown in [16] that by constructing the epidemic and random graph in unison, the generalized stochastic epidemic upon a Bernoulli random graph also satisfies the randomized Reed–Frost criterion.

---

---

The asymptotic and exact final size distribution of the randomized Reed–Frost epidemic have been derived in [14] and [18], respectively. The aim of the current work is to obtain the asymptotic final size distribution of a *multitype randomized Reed–Frost* epidemic. (The exact final size distribution is given in [18].) A full description of the multitype randomized Reed–Frost epidemic will be given in Section 2. Throughout we shall take the final size of the epidemic to denote the total number of initially susceptible individuals infected during the course of the epidemic. The analysis of the asymptotic final size distribution involves a method of proof which differs markedly from that given in [14]. We utilize an embedding argument based upon [19] and [21]; see also [6, 8, 9, 10]. In particular, our results generalize those of [6], Section 4. Conventionally (see, e.g., [6]), an individual's type is determined deterministically. However, we also consider the case where the allocation of individuals to type is random. Thus the model considered in [7] is another special case.

The main motivation for the current work was extending the asymptotic results of [16] to epidemics upon more general random graphs/population networks. Epidemics upon graphs have received considerable attention in recent years; see, for example, [3, 4, 16, 17]. Only in [16] is a central limit theorem for the final size of the epidemic considered. Thus in Section 7.3, we show that epidemics upon various extensions of the Bernoulli random graph satisfy the multitype randomized Reed–Frost criterion, and therefore we obtain results which are considerably more extensive than those in [16].

The paper is structured as follows. In Section 2 the multitype randomized Reed–Frost epidemic is defined, and a construction of the epidemic process suitable for analysis is given. In Section 3 a weak law of large numbers result is obtained for the proportion of the initially susceptible population that are ultimately infected by the epidemic. This is followed in Section 4 by a branching process approximation to the epidemic. This is appropriate when the total number of initial infectives is small. The main result is presented in Section 5, and is the derivation of an asymptotic Gaussian limit for the final size of the epidemic. While the results proved in Sections 3–5 are "standard" for epidemic processes, the method of proof, especially in Section 5, is novel. We present a simpler proof for deriving the asymptotic Gaussian limit for the final size distribution. This proof is applicable for a wide range of epidemic models. In particular, the approach taken in Section 5 can be applied to the great circle model of [10] to extend the results therein. Throughout Sections 3–5, it will be assumed that the allocation of individuals to types is deterministic. Therefore in Section 6 we outline the (minor) modifications required to consider random allocation of individuals to types. Finally in Section 7 the results are specialized to the models of [6] and [7], and epidemics upon extensions of the Bernoulli random graph.

**2. Description and construction of the epidemic process.** Consider a closed population divided into $m$ types (or groups). Let $\mathbf{a} = (a_1, a_2, \ldots, a_m)$ and $\mathbf{N} = (N_1, N_2, \ldots, N_m)$ denote the total number of initial infectives and susceptibles, respectively. [We shall take the statement that a group of individuals, $S$ say, is of size $\mathbf{s} = (s_1, s_2, \ldots, s_m)$, to be interpreted that $S$ contains $s_k$ individuals of type $k$ $(k = 1, 2, \ldots, m)$.] For $i = 1, 2, \ldots, m$, let $\mathbf{V}_i = (V_{i,1}, V_{i,2}, \ldots, V_{i,m})$ be an $m$-dimensional random vector taking values on the $m$-dimensional unit cube $[\mathbf{0}, \mathbf{1}]$, where $\mathbf{0}$ and $\mathbf{1}$ denote the $m$-dimensional vectors all of whose components are 0 and 1, respectively. [Throughout this paper, all vectors are row vectors, and inequalities between vectors are to be interpreted componentwise, with $\mathbf{u} \leq \mathbf{v}$ ($\mathbf{u} < \mathbf{v}$), if for $1 \leq i \leq m$, $u_i \leq v_i$ (and at least one of the inequalities is strict). Furthermore, for any vectors $\mathbf{u}, \mathbf{v}, \mathbf{w} \in \mathbb{R}^m$, $[\mathbf{u}, \mathbf{v}] = \{\mathbf{w} : \mathbf{u} \leq \mathbf{w} \leq \mathbf{v}\}$.] The infectious behaviors of all individuals are assumed to be independent. For $i = 1, 2, \ldots, m$, assign to each infective of type-$i$ an independent infectivity random vector distributed according to $\mathbf{V}_i$. Then for $i = 1, 2, \ldots, m$, a type-$i$ infective with infectivity random vector $\mathbf{V}_{(i,1)} = (V_{(i,1),1}, V_{(i,1),2}, \ldots, V_{(i,1),m})$ has probability $V_{(i,1),k}$ of making, while infectious, an infectious contact with a given susceptible of type $k$, independent of whether the individual makes any other infectious contacts or not. This gives a natural multitype extension of the randomized Reed–Frost epidemic described on page 270 of [14]. Clearly, the total number of individuals of each type contacted by a given infective can be dependent. This is in contrast to [2], where a multitype version of the Martin–Löf [14] model is considered.

For the asymptotic results derived in this paper, it will be convenient to consider a sequence of epidemics $\{E_\nu\}$ as $\nu \to \infty$. Consider fixed $\nu$. We shall assume that $\mathbf{a}^\nu$, $\mathbf{N}^\nu$ and $\mathbf{V}_i^\nu$ $(1 \leq i \leq m)$ are dependent upon $\nu$. Let $N^\nu = N_1^\nu + N_2^\nu + \cdots + N_m^\nu$, with $N^\nu \to \infty$ as $\nu \to \infty$. For $i = 1, 2, \ldots, m$, let $\pi_i^\nu = \frac{N_i^\nu}{N^\nu} \to \pi_i > 0$, and $\zeta_i^\nu = \frac{a_i^\nu}{N^\nu \pi_i} \to \zeta_i \geq 0$ as $\nu \to \infty$. For $i = 1, 2, \ldots, m$, label the initial infectives of type $i$, individuals $(i, -(a_i^\nu - 1)), (i, -(a_i^\nu - 2)), \ldots, (i, 0)$, and label the initial susceptibles of type $i$, individuals $(i, 1), (i, 2), \ldots, (i, N_i^\nu)$. Assign to the initial infectives of type $i$, infectivity random vectors $\mathbf{V}_{(i,1)}^\nu$, $\mathbf{V}_{(i,2)}^\nu, \ldots, \mathbf{V}_{(i,a_i^\nu)}^\nu$, and for initially susceptible individuals of type $i$ infected during the course of the epidemic, assign infectivity random vectors sequentially from $\mathbf{V}_{(i,a_i^\nu+1)}^\nu, \mathbf{V}_{(i,a_i^\nu+2)}^\nu, \ldots, \mathbf{V}_{(i,a_i^\nu+N_i^\nu)}^\nu$.

For each $\nu \geq 1$, the epidemic process can be constructed on a generational basis by sampling suitably defined stopping times of a simple multivariate counting process, $\mathbf{X}^\nu(\cdot)$. The embedding of the epidemic processes (and in particular, its final size distribution) into the $\mathbf{X}^\nu(\cdot)$-processes is based on the approach introduced in [19]; see also [6, 9, 21].

For $\nu \geq 1$ and $\mathbf{t} \geq 0$, consider the first $[\mathbf{t}N^\nu\Pi]$ infectives with infectivity random vectors,

$$\{\mathbf{V}^\nu_{(1,1)}, \mathbf{V}^\nu_{(1,2)}, \ldots, \mathbf{V}^\nu_{(1,[t_1N^\nu\pi_1])}, \mathbf{V}^\nu_{(2,1)}, \mathbf{V}^\nu_{(2,2)}, \ldots, \mathbf{V}^\nu_{(m,[t_mN^\nu\pi_m])}\},$$

where $\Pi = \operatorname{diag}(\pi_1, \pi_2, \ldots, \pi_m)$. Then for $1 \leq i \leq m$, and $1 \leq j \leq N^\nu_i$, let $\chi^\nu_{i,j}(\mathbf{t}) = 1$, if individual $(i,j)$ is infected when exposed to the first $[\mathbf{t}N^\nu\Pi]$ infectives, and $\chi^\nu_{i,j}(\mathbf{t}) = 0$, otherwise. Hence for fixed $1 \leq i \leq m$ and $\mathbf{t} \geq \mathbf{0}$, the $\{\chi^\nu_{i,j}(\mathbf{t})\}$ are identically distributed. For $i = 1, 2, \ldots, m$, and $\mathbf{t} \geq \mathbf{0}$, let $X^\nu_i(\mathbf{t}) = \sum_{j=1}^{N^\nu_i} \chi^\nu_{i,j}(\mathbf{t})$, with $\mathbf{X}^\nu(\mathbf{t}) = (X^\nu_1(\mathbf{t}), X^\nu_2(\mathbf{t}), \ldots, X^\nu_m(\mathbf{t}))$. Thus $\mathbf{X}^\nu(\mathbf{t})$ counts the total number of initial susceptibles that would be infected if each individual were exposed to the first $[\mathbf{t}N^\nu\Pi]$ infectives.

The epidemic process $E^\nu$ can be constructed using the $\mathbf{X}^\nu(\cdot)$-processes as follows. Let generation 0 denote the initial infectives, and for $k \geq 1$, let generation $k$ denote those individuals infected by the infectives in generation $k - 1$. For $k \geq 1$, let $\mathbf{T}^\nu_k$ denote the total number of initial susceptibles infected by the infectives in generations 0 through to $k - 1$, inclusive. Therefore for $k \geq 0$, $\mathbf{T}^\nu_{k+1} = \mathbf{X}^\nu(\bar{\mathbf{T}}^\nu_k + \boldsymbol{\zeta}^\nu_k)$, where $\mathbf{T}^\nu_0 = \mathbf{0}$, and for $j \geq 0$, $\bar{\mathbf{T}}^\nu_j = \mathbf{T}^\nu_j(N^\nu\Pi)^{-1}$. Thus for $k \geq 1$,

$$\bar{\mathbf{T}}^\nu_{k+1} = \mathbf{X}^\nu(\bar{\mathbf{T}}^\nu_k + \boldsymbol{\zeta}^\nu_k)(N^\nu\Pi)^{-1}.$$

Consequently, the above sequence stops at generation $k^*$, where $k^* = \min\{k : \bar{\mathbf{T}}^\nu_{k+1} = \bar{\mathbf{T}}^\nu_k\}$. Since $N^\nu$ is finite, $k^*$ is well defined. Let $\mathbf{T}^\nu_\infty(= \bar{\mathbf{T}}^\nu_\infty(N^\nu\Pi))$ denote the final size of the epidemic. Then $\bar{\mathbf{T}}^\nu_\infty = \bar{\mathbf{T}}^\nu_{k^*}$, and satisfies

$$(2.1) \qquad \bar{\mathbf{T}}^\nu_\infty = \mathbf{X}^\nu(\bar{\mathbf{T}}^\nu_\infty + \boldsymbol{\zeta}^\nu)(N^\nu\Pi)^{-1}.$$

Finally in Section 6 we shall consider the case where for each $\nu \geq 1$, $\mathbf{a}^\nu = N^\nu\boldsymbol{\zeta}^\nu$ and $\mathbf{N}^\nu = N^\nu\boldsymbol{\pi}^\nu$ are random vectors, with $a^\nu$ and $N^\nu$ constant. The above construction of the epidemic process can again be used with the assumption that $\boldsymbol{\zeta}^\nu \xrightarrow{p} \boldsymbol{\zeta}$ and $\boldsymbol{\pi}^\nu \xrightarrow{p} \boldsymbol{\pi}$ as $\nu \to \infty$.

**3. Weak law of large numbers for the final size distribution.** In this section we shall consider the limiting distribution of $\bar{\mathbf{T}}^\nu_\infty$ as $\nu \to \infty$. In particular, we show in Corollary 3.3 that

$$\min\{|\bar{\mathbf{T}}^\nu_\infty|, |\bar{\mathbf{T}}^\nu_\infty - \boldsymbol{\tau}|\} \xrightarrow{p} 0 \qquad \text{as } \nu \to \infty,$$

where $\boldsymbol{\tau}$ is the solution of an appropriate deterministic model.

We shall begin by deriving $\boldsymbol{\tau}$. For $1 \leq i, j, k \leq m$, let $\mathbb{E}[V^\nu_{i,k}] = \mu^\nu_{ik}$, and let $\operatorname{cov}(V^\nu_{i,j}, V^\nu_{i,k}) = \lambda^\nu_{ijk}$. We shall require that for $1 \leq i, k \leq m$, there exists $0 \leq \mu_{ik} < \infty$, such that $N^\nu\mu^\nu_{ik} \to \mu_{ik}$ and $N^\nu\lambda^\nu_{ikk} \to 0$ as $\nu \to \infty$.

For $\nu \geq 1$, let

$$\tilde{\mathbf{V}}^\nu = \{\mathbf{V}^\nu_{(1,1)}, \mathbf{V}^\nu_{(1,2)}, \ldots, \mathbf{V}^\nu_{(1,a^\nu_1+N^\nu_1)}, \mathbf{V}^\nu_{(2,1)}, \mathbf{V}^\nu_{(2,2)}, \ldots, \mathbf{V}^\nu_{(m,a^\nu_m+N^\nu_m)}\}$$

denote the infectivity random vectors of all individuals in the population. For $\nu \geq 1$ and $1 \leq i \leq m$, let $r_i^\nu(\mathbf{t}) = \mathbb{E}[\chi_{i,1}^\nu(\mathbf{t})]$ and $r_i^\nu(\mathbf{t}|\tilde{\mathbf{V}}^\nu) = \mathbb{E}[\chi_{i,1}^\nu(\mathbf{t})|\tilde{\mathbf{V}}^\nu]$. Let $\tilde{N}_i^\nu = N^\nu \pi_i$ $(1 \leq i \leq m)$. Then for $1 \leq i \leq m$,

$$r_i^\nu(\mathbf{t}|\tilde{\mathbf{V}}^\nu) = \mathbb{E}[\chi_{i,1}^\nu(\mathbf{t})|\tilde{\mathbf{V}}^\nu] = 1 - \prod_{k=1}^m \prod_{l=1}^{[t_k \tilde{N}_k^\nu]} (1 - V_{(k,l),i}^\nu)$$

and

$$r_i^\nu(\mathbf{t}) = \mathbb{E}_{\tilde{\mathbf{V}}^\nu}[r_i^\nu(\mathbf{t}|\tilde{\mathbf{V}}^\nu)] = 1 - \prod_{k=1}^m (1 - \mu_{ki}^\nu)^{[t_k \tilde{N}_k^\nu]}.$$

For $1 \leq i \leq m$ and $\mathbf{t} \geq \mathbf{0}$, let $r_i(\mathbf{t}) = 1 - \exp(-\sum_{k=1}^m t_k \pi_k \mu_{ki})$. It is trivial to show that for $1 \leq i \leq m$ and $\mathbf{t} \geq \mathbf{0}$, $r_i^\nu(\mathbf{t}) \to r_i(\mathbf{t})$ as $\nu \to \infty$.

For $\boldsymbol{\zeta} \geq \mathbf{0}$, let

$$\mathcal{T}(\boldsymbol{\zeta}) = \{\boldsymbol{\tau} \geq \mathbf{0} : \boldsymbol{\tau} = \mathbf{r}(\boldsymbol{\tau} + \boldsymbol{\zeta})\},$$

where $\mathbf{r}(\mathbf{t}) = (r_1(\mathbf{t}), r_2(\mathbf{t}), \ldots, r_m(\mathbf{t}))$ $(\mathbf{t} \geq 0)$. Let $\boldsymbol{\sigma} = \mathbf{1} - \boldsymbol{\tau}$. Then $\boldsymbol{\sigma} = (\sigma_1, \sigma_2, \ldots, \sigma_m)$ satisfies

$$(3.1) \qquad \sigma_i = \exp\left(-\sum_{k=1}^m \pi_k \mu_{ki}(1 + \zeta_k - \sigma_k)\right) \qquad (i = 1, 2, \ldots, m).$$

Equation (3.1) is identical to [20], equation (4) and [6], equation (2.2). Therefore denoting the matrix with elements $\mu_{ki}$ by $M$, it follows from [20], Lemma 1, that if $M^T \Pi$ is irreducible, and $\boldsymbol{\zeta} > \mathbf{0}$, then the equation (3.1) for $\boldsymbol{\sigma}$, and hence, $\boldsymbol{\tau}$ has a unique solution in $[\mathbf{0}, \mathbf{1}]$. If $\boldsymbol{\zeta} = \mathbf{0}$, then $\boldsymbol{\tau} = \mathbf{0}$ is a solution of (3.1). Let $R$ denote the Perron–Frobenius eigenvalue of $M\Pi$. Then for $\boldsymbol{\zeta} = \mathbf{0}$, there exists $\boldsymbol{\tau} \in \mathcal{T}(\mathbf{0})$ such that $\boldsymbol{\tau} > \mathbf{0}$, if and only if $R > 1$. We shall let $\boldsymbol{\tau}$ denote the nonzero solution of $\boldsymbol{\tau} = \mathbf{r}(\boldsymbol{\tau} + \boldsymbol{\zeta})$ whenever such a solution exists (i.e., $\boldsymbol{\zeta} > \mathbf{0}$ or $R > 1$).

Before analyzing the $\mathbf{X}^\nu(\cdot)$-processes (Lemma 3.2), and hence $\bar{\mathbf{T}}_\infty^\nu$ (Corollary 3.3), we develop two useful bounds for the covariances of $\{\chi_{i,j}^\nu(\cdot)\}$ $(1 \leq i \leq m)$. The bound (3.3) will not be required until Section 5.

LEMMA 3.1. *For $1 \leq i \leq m$, $1 \leq j < l \leq N_i^\nu$ and for all $\mathbf{t}, \mathbf{u} \geq \mathbf{0}$,*

$$(3.2) \qquad 0 \leq \mathrm{cov}(\chi_{i,j}^\nu(\mathbf{t}), \chi_{i,l}^\nu(\mathbf{u})) \leq \sum_{k=1}^m (t_k \wedge u_k) \tilde{N}_k^\nu \lambda_{kii}^\nu$$

*and*

$$|\mathrm{cov}(\chi_{i,j}^\nu(\mathbf{t}) - \chi_{i,j}^\nu(\mathbf{u}), \chi_{i,l}^\nu(\mathbf{t}) - \chi_{i,l}^\nu(\mathbf{u}))|$$
$$(3.3) \qquad \leq 6\left\{\sum_{k=1}^m ([(t_k \vee u_k)\tilde{N}_k^\nu] - [(t_k \wedge u_k)\tilde{N}_k^\nu])\mu_{ki}^\nu\right\}\left\{\sum_{k=1}^m (t_k \wedge u_k)\tilde{N}_k^\nu \lambda_{kii}^\nu\right\}$$

$$+ 2 \sum_{k=1}^{m} ([(t_k \vee u_k) \tilde{N}_k^{\nu}] - [(t_k \wedge u_k) \tilde{N}_k^{\nu}]) \lambda_{kii}^{\nu}.$$

PROOF. Fix $\mathbf{t}, \mathbf{u} \geq \mathbf{0}$. Let $\mathbf{v} = (t_1 \wedge u_1, t_2 \wedge u_2, \dots, t_m \wedge u_m)$ and $\mathbf{w} = (t_1 \vee u_1, t_2 \vee u_2, \dots, t_m \vee u_m)$. For $1 \leq i \leq m$ and $1 \leq j < l \leq N_i^{\nu}$,

$$\operatorname{cov}(\chi_{i,j}^{\nu}(\mathbf{t}), \chi_{i,l}^{\nu}(\mathbf{u})) = \operatorname{cov}(1 - \chi_{i,j}^{\nu}(\mathbf{t}), 1 - \chi_{i,l}^{\nu}(\mathbf{u}))$$

with

$$\mathbb{E}[\{1 - \chi_{i,j}^{\nu}(\mathbf{t})\}\{1 - \chi_{i,l}^{\nu}(\mathbf{u})\}]$$
$$= \prod_{k=1}^{m} \left\{ \prod_{l=1}^{[v_k \tilde{N}_k^{\nu}]} \mathbb{E}[(1 - V_{(k,l),i}^{\nu})^2] \prod_{l=[v_k \tilde{N}_k^{\nu}]+1}^{[w_k \tilde{N}_k^{\nu}]} \mathbb{E}[1 - V_{(k,l),i}^{\nu}] \right\}$$

and $\prod_{l=a}^{b}(\cdot) = 1$ if $a > b$.

Therefore after some straightforward algebra, we have that

$$\begin{aligned}
(3.4) \quad & \operatorname{cov}(\chi_{i,j}^{\nu}(\mathbf{t}), \chi_{i,l}^{\nu}(\mathbf{u})) \\
&= \left\{ \prod_{k=1}^{m}(1 - \mu_{ki}^{\nu})^{[w_k \tilde{N}_k^{\nu}] - [v_k \tilde{N}_k^{\nu}]} \right\} \\
&\quad \times \left( \prod_{k=1}^{m} \{(1 - \mu_{ki}^{\nu})^2 + \lambda_{kii}^{\nu}\}^{[v_k \tilde{N}_k^{\nu}]} - \prod_{k=1}^{m} \{(1 - \mu_{ki}^{\nu})^2\}^{[v_k \tilde{N}_k^{\nu}]} \right).
\end{aligned}$$

It follows trivially from (3.4) that $\operatorname{cov}(\chi_{i,j}^{\nu}(\mathbf{t}), \chi_{i,l}^{\nu}(\mathbf{u})) \geq 0$. Since all the products in (3.4) are less than or equal to 1, we have that

$$\begin{aligned}
(3.5) \quad & \operatorname{cov}(\chi_{i,j}^{\nu}(\mathbf{t}), \chi_{i,l}^{\nu}(\mathbf{u})) \\
&\leq \sum_{k=1}^{m}[v_k \tilde{N}_k^{\nu}]\{(1 - \mu_{ki}^{\nu})^2 + \lambda_{kii}^{\nu} - (1 - \mu_{ki}^{\nu})^2\} = \sum_{k=1}^{m}[v_k \tilde{N}_k^{\nu}]\lambda_{kii}^{\nu}.
\end{aligned}$$

Hence (3.2) is proved.

Straightforward algebraic manipulation of (3.4) gives

$$\begin{aligned}
& \operatorname{cov}(\chi_{i,j}^{\nu}(\mathbf{t}), \chi_{i,l}^{\nu}(\mathbf{t})) - \operatorname{cov}(\chi_{i,j}^{\nu}(\mathbf{t}), \chi_{i,l}^{\nu}(\mathbf{u})) \\
&= \left( \prod_{k=1}^{m} \{(1 - \mu_{ki}^{\nu})^2 + \lambda_{kii}^{\nu}\}^{[t_k \tilde{N}_k^{\nu}] - [v_k \tilde{N}_k^{\nu}]} - \prod_{k=1}^{m}(1 - \mu_{ki}^{\nu})^{[w_k \tilde{N}_k^{\nu}] - [v_k \tilde{N}_k^{\nu}]} \right) \\
(3.6) \quad & \quad \times \operatorname{cov}(\chi_{i,j}^{\nu}(\mathbf{v}), \chi_{i,l}^{\nu}(\mathbf{v})) \\
& \quad + \left\{ \prod_{k=1}^{m}(1 - \mu_{ki}^{\nu})^{2[v_k \tilde{N}_k^{\nu}]} \right\}
\end{aligned}$$

$$\times \left( \prod_{k=1}^{m} \{(1-\mu_{ki}^{\nu})^2 + \lambda_{kii}^{\nu}\}^{[t_k \tilde{N}_k^{\nu}] - [v_k \tilde{N}_k^{\nu}]} \right.$$

$$\left. - \prod_{k=1}^{m} \{(1-\mu_{ki}^{\nu})^2\}^{[t_k \tilde{N}_k^{\nu}] - [v_k \tilde{N}_k^{\nu}]} \right).$$

By arguments identical to those used in (3.5),

$$(3.7) \quad \begin{aligned} &\left| \prod_{k=1}^{m} \{(1-\mu_{ki}^{\nu})^2 + \lambda_{kii}^{\nu}\}^{[t_k \tilde{N}_k^{\nu}] - [v_k \tilde{N}_k^{\nu}]} - \prod_{k=1}^{m} (1-\mu_{ki}^{\nu})^{[w_k \tilde{N}_k^{\nu}] - [v_k \tilde{N}_k^{\nu}]} \right| \\ &\leq \left| \prod_{k=1}^{m} \{(1-\mu_{ki}^{\nu})^2\}^{[t_k \tilde{N}_k^{\nu}] - [v_k \tilde{N}_k^{\nu}]} - 1 \right| \\ &\quad + \left| 1 - \prod_{k=1}^{m} (1-\mu_{ki}^{\nu})^{[w_k \tilde{N}_k^{\nu}] - [v_k \tilde{N}_k^{\nu}]} \right| \\ &\leq 2 \sum_{k=1}^{m} \{[t_k \tilde{N}_k^{\nu}] - [v_k \tilde{N}_k^{\nu}]\} \mu_{ki}^{\nu} + \sum_{k=1}^{m} \{[w_k \tilde{N}_k^{\nu}] - [v_k \tilde{N}_k^{\nu}]\} \mu_{ki}^{\nu}. \end{aligned}$$

A similar inequality exists for the latter term on the right-hand side of (3.6) giving

$$(3.8) \quad \begin{aligned} &|\operatorname{cov}(\chi_{i,j}^{\nu}(\mathbf{t}), \chi_{i,l}^{\nu}(\mathbf{t})) - \operatorname{cov}(\chi_{i,j}^{\nu}(\mathbf{t}), \chi_{i,l}^{\nu}(\mathbf{u}))| \\ &\leq 3\left\{ \sum_{k=1}^{m} ([(t_k \vee u_k) \tilde{N}_k^{\nu}] - [(t_k \wedge u_k) \tilde{N}_k^{\nu}]) \mu_{ki}^{\nu} \right\} \left\{ \sum_{k=1}^{m} (t_k \wedge u_k) \tilde{N}_k^{\nu} \lambda_{kii}^{\nu} \right\} \\ &\quad + \sum_{k=1}^{m} ([(t_k \vee u_k) \tilde{N}_k^{\nu}] - [(t_k \wedge u_k) \tilde{N}_k^{\nu}]) \lambda_{kii}^{\nu}. \end{aligned}$$

Therefore, since the inequality in (3.8) holds with $\mathbf{t}$ and $\mathbf{u}$ reversed, (3.3) follows immediately. $\quad\square$

LEMMA 3.2. *For any $\mathbf{s} \geq \mathbf{0}$, let $D(\mathbf{s}) = [\mathbf{0}, \mathbf{s}]$. For $1 \leq i \leq m$, and for any $\mathbf{s} \geq \mathbf{0}$,*

$$(3.9) \qquad \sup_{\mathbf{t} \in D(\mathbf{s})} \left| \frac{1}{\tilde{N}_i^{\nu}} X_i^{\nu}(\mathbf{t}) - r_i(\mathbf{t}) \right| \xrightarrow{p} 0 \qquad \text{as } \nu \to \infty.$$

*Thus*

$$(3.10) \qquad \sup_{\mathbf{t} \in D(\mathbf{s})} |\mathbf{X}^{\nu}(\mathbf{t})(N^{\nu}\Pi)^{-1} - \mathbf{r}(\mathbf{t})| \xrightarrow{p} 0 \qquad \text{as } \nu \to \infty.$$

Proof. Consider fixed $1 \leq i \leq m$ and $\mathbf{t} \in D(\mathbf{s})$. Since $r_i^\nu(\mathbf{t}) \to r_i(\mathbf{t})$ and $\pi_i^\nu \to \pi$ as $\nu \to \infty$, to prove that $\frac{1}{N_i^\nu} X_i^\nu(\mathbf{t}) \xrightarrow{p} r_i(\mathbf{t})$, it is sufficient to show that

$$\frac{1}{N_i^\nu} X_i^\nu(\mathbf{t}) - \frac{\pi_i^\nu}{\pi_i^\nu} r_i^\nu(\mathbf{t}) = \frac{1}{N_i^\nu} X_i^\nu(\mathbf{t}) - \frac{N_i^\nu}{N_i^\nu} r_i^\nu(\mathbf{t}) \xrightarrow{p} 0 \qquad \text{as } \nu \to \infty.$$

Fix $\varepsilon > 0$. By Chebyshev's inequality,

$$\mathbb{P}\left( \left| \frac{1}{\tilde{N}_i^\nu} X_i^\nu(\mathbf{t}) - \frac{N_i^\nu}{\tilde{N}_i^\nu} r_i^\nu(\mathbf{t}) \right| > \varepsilon \right) \leq \frac{1}{\varepsilon^2 (\tilde{N}_i^\nu)^2} \sum_{j=1}^{N_i^\nu} \sum_{l=1}^{N_i^\nu} \operatorname{cov}(\chi_{i,j}^\nu(\mathbf{t}), \chi_{i,l}^\nu(\mathbf{t})).$$

Therefore, by (3.2),

$$\mathbb{P}\left( \left| \frac{1}{N_i^\nu} X_i^\nu(\mathbf{t}) - \frac{N_i^\nu}{\tilde{N}_i^\nu} r_i^\nu(\mathbf{t}) \right| > \varepsilon \right)$$

$$\leq \frac{1}{\varepsilon^2} \left\{ \frac{N_i^\nu}{(\tilde{N}_i^\nu)^2} r_i^\nu(\mathbf{t})(1 - r_i^\nu(\mathbf{t})) + \frac{N_i^\nu(N_i^\nu - 1)}{(\tilde{N}_i^\nu)^2} \sum_{k=1}^{m} t_k \tilde{N}_k^\nu \lambda_{kii}^\nu \right\},$$

and so, for all $\mathbf{t} \in D(\mathbf{s})$,

$$(3.11) \qquad \left| \frac{1}{\tilde{N}_i^\nu} X_i^\nu(\mathbf{t}) - r_i(\mathbf{t}) \right| \xrightarrow{p} 0 \qquad \text{as } \nu \to \infty.$$

The remainder of the proof of (3.9) follows by arguments similar to [8], Lemma 5.1 since both $X_i^\nu(\mathbf{t})$ and $r_i(\mathbf{t})$ are nondecreasing in $\mathbf{t}$. Therefore the details are omitted. □

Corollary 3.3 follows immediately from (2.1), Lemma 3.2 and (3.1).

Corollary 3.3. *Suppose that for $1 \leq i \leq m$, $\boldsymbol{\pi}^\nu \to \boldsymbol{\pi}$ and $\boldsymbol{\zeta}^\nu \to \boldsymbol{\zeta}$ as $\nu \to \infty$. Then if $\boldsymbol{\zeta} > \mathbf{0}$,*

$$\bar{\mathbf{T}}_\infty^\nu \xrightarrow{p} \boldsymbol{\tau} \qquad \text{as } \nu \to \infty.$$

*Alternatively, if $\boldsymbol{\zeta} = \mathbf{0}$,*

$$(3.12) \qquad \min\{|\bar{\mathbf{T}}_\infty^\nu|, |\bar{\mathbf{T}}_\infty^\nu - \boldsymbol{\tau}|\} \xrightarrow{p} 0 \qquad \text{as } \nu \to \infty,$$

*where $\boldsymbol{\tau}$ exists if and only if $R > 1$.*

**4. Branching process approximation.** We shall assume that there exists $\mathbf{a} > 0$ such that for all $\nu \geq 1$, $\mathbf{a}^\nu = \mathbf{a}$. (The extension to the case $\mathbf{a}^\nu \to \mathbf{a}$ as $\nu \to \infty$ is trivial.) Corollary 3.3, (3.12) is valid in this situation and therefore our aim is twofold: first, to find $\mathbb{P}(\lim_{\nu \to \infty} \bar{\mathbf{T}}_\infty^\nu = \boldsymbol{\tau})$, that is, the probability of a major epidemic outbreak, and second, to find the limiting distribution of $\mathbf{T}_\infty^\nu$ conditional upon a small epidemic outbreak occurring (i.e., $\bar{\mathbf{T}}_\infty^\nu \xrightarrow{p} 0$ as $\nu \to \infty$).

We shall show that $\mathbf{T}_\infty^\nu$ converges in distribution to the total progeny $\mathbf{Z}$ of a suitably defined multitype Galton–Watson branching process. We proceed by giving a slightly different construction of the epidemic process from that given in Section 2. It is straightforward to show that the two constructions are equivalent, with regard to the final size of the epidemic process. We shall then describe the branching processes to which the epidemic processes are coupled, along with the limiting branching process which has total progeny $\mathbf{Z}$.

We make the following assumption. For $1 \le k \le m$, there exists a random variable $\mathbf{U}_k$ such that

$$\mathbf{V}_k^\nu N^\nu \Pi^\nu \xrightarrow{D} \mathbf{U}_k = (U_{k,1}, U_{k,2}, \ldots, U_{k,m}) \qquad \text{as } \nu \to \infty,$$

with $\mathbb{E}[U_{k,i}] = \mu_{ki} < \infty$ $(1 \le i \le m)$. For $1 \le k \le m$, let $\mathbf{R}_k^\nu (= (R_{k,1}^\nu, R_{k,2}^\nu, \ldots, R_{k,m}^\nu)) = (\text{Bin}(N_1^\nu, V_{k,1}^\nu), \text{Bin}(N_2^\nu, V_{k,2}^\nu), \ldots, \text{Bin}(N_m^\nu, V_{k,m}^\nu))$, with $\tilde{\mathbf{R}}^\nu = (\mathbf{R}_1^\nu, \mathbf{R}_2^\nu, \ldots, \mathbf{R}_m^\nu)$. Thus for $1 \le i, k \le m$, let $R_{k,i}^\nu$ denote the total number of initially susceptible individuals of type $i$ contacted by a typical type-$k$ infective. For $1 \le i, k \le m$ and $l \ge 1$, let $P_{(k,l),i}^\nu$ denote the set of type-$i$ individuals contacted by the $l$th infective of type $k$, with $R_{(k,l),i}^\nu = |P_{(k,l),i}^\nu|$. Finally, let $\mathbf{P}_{(k,l)}^\nu = \{P_{(k,l),1}^\nu, P_{(k,l),2}^\nu, \ldots, P_{(k,l),m}^\nu\}$.

Let $\mathbf{Z}^\nu = (Z_1^\nu, Z_2^\nu, \ldots, Z_m^\nu)$ denote the total progeny, excluding the $\mathbf{a}$ initial ancestors, of an $m$-type Galton–Watson branching process with (mixed binomial) reproductive law $\tilde{\mathbf{R}}^\nu$. That is, a typical type-$k$ individual has $R_{k,i}^\nu$ offspring of type $i$. For $1 \le k \le m$, let $\mathbf{R}_k = (Po(U_{k,1}), Po(U_{k,2}), \ldots, Po(U_{k,m}))$. Therefore let $\mathbf{Z} = (Z_1, Z_2, \ldots, Z_m)$ denote the total progeny, excluding the $\mathbf{a}$ initial ancestors, of an $m$-type Galton–Watson branching process with (mixed Poisson) reproductive law, $\tilde{\mathbf{R}} = (\mathbf{R}_1, \mathbf{R}_2, \ldots, \mathbf{R}_m)$.

LEMMA 4.1. *Suppose that for $1 \le k \le m$, $\mathbf{V}_k^\nu N^\nu \Pi^\nu \xrightarrow{D} \mathbf{U}_k$ as $\nu \to \infty$. Then:*

(i) $\mathbf{R}_k^\nu \xrightarrow{D} \mathbf{R}_k$ *as $\nu \to \infty$,*

(ii) $\mathbf{T}_\infty^\nu \xrightarrow{D} \mathbf{Z}$ *as $\nu \to \infty$.*

PROOF. (i) The result is a trivial extension of [13], Lemma 5.8(i). Hence the details are omitted.

(ii) Note that an immediate consequence of (i) is that $\mathbf{Z}^\nu \xrightarrow{D} \mathbf{Z}$ as $\nu \to \infty$. For $\mathbf{0} \le \mathbf{s} \in \mathbb{Z}^m$, let $M^\nu(\mathbf{s})$ denote the event that the infectives

$$(1,1), (1,2), \ldots, (1, s_1), (2,1), (2,2), \ldots, (2, s_2), \ldots, (m, s_m)$$

make infectious contacts with distinct sets of individuals. The key point is that conditional upon $M^\nu(\mathbf{s})$ occurring, the epidemic and branching processes can be coupled, so that the processes coincide for the first $\mathbf{s}$ individuals. Now for $\mathbf{s} \geq \mathbf{0}$, it is straightforward to show that

$$\mathbb{P}(M^\nu(\mathbf{s})^C) = \mathbb{P}\left(\left\{\bigcup_{i=1}^{m-1}\bigcup_{k=i+1}^{m}\bigcup_{j=1}^{s_i}\bigcup_{l=1}^{s_k}\{\mathbf{P}^\nu_{(i,j)} \cap \mathbf{P}^\nu_{(k,l)}\}\right\}\right.$$

$$\left.\cup\left\{\bigcup_{i=1}^{m}\bigcup_{j=1}^{s_i-1}\bigcup_{l=j+1}^{s_i}\{\mathbf{P}^\nu_{(i,j)} \cap \mathbf{P}^\nu_{(i,l)}\}\right\} \neq \varnothing\right)$$

(4.1)

$$\leq \sum_{i=1}^{m-1}\sum_{k=i+1}^{m}\sum_{j=1}^{s_i}\sum_{l=1}^{s_k}\mathbb{P}(\mathbf{P}^\nu_{(i,j)} \cap \mathbf{P}^\nu_{(k,l)} \neq \varnothing)$$

$$+ \sum_{i=1}^{m}\sum_{j=1}^{s_i-1}\sum_{l=j+1}^{s_i}\mathbb{P}(\mathbf{P}^\nu_{(i,j)} \cap \mathbf{P}^\nu_{(i,l)} \neq \varnothing).$$

For $(i,j) \neq (k,l)$, $1 \leq d \leq m$ and $b,c \geq 0$,

$$\mathbb{P}(P^\nu_{(i,j),d} \cap P^\nu_{(k,l),d} \neq \varnothing | R^\nu_{(i,j),d} = b, R^\nu_{(k,l),d} = c) \leq \frac{bc}{N^\nu_d}.$$

Thus for $(i,j) \neq (k,l)$,

$$\mathbb{P}(\mathbf{P}^\nu_{(i,j)} \cap \mathbf{P}^\nu_{(k,l)} \neq \varnothing) \leq \sum_{d=1}^{m}\frac{1}{N^\nu_d}\mathbb{E}[R^\nu_{(i,j),d}]\mathbb{E}[R^\nu_{(k,l),d}] \to 0 \qquad \text{as } \nu \to \infty.$$

Therefore the right-hand side of (4.1) converges to 0 as $\nu \to \infty$.

Thus for all $\mathbf{s} \in (\mathbb{Z}^+)^m$,

$$\begin{aligned}\mathbb{P}(\mathbf{T}^\nu_\infty \leq \mathbf{s}) &= \mathbb{P}(\mathbf{T}^\nu_\infty \leq \mathbf{s}|M^\nu(\mathbf{s}+\mathbf{a}))\mathbb{P}(M^\nu(\mathbf{s}+\mathbf{a}))\\&\quad + \mathbb{P}(\mathbf{T}^\nu_\infty \leq \mathbf{s}|M^\nu(\mathbf{s}+\mathbf{a})^C)\mathbb{P}(M^\nu(\mathbf{s}+\mathbf{a})^C)\\&= \mathbb{P}(\mathbf{Z}^\nu \leq \mathbf{s}|M^\nu(\mathbf{s}+\mathbf{a}))\mathbb{P}(M^\nu(\mathbf{s}+\mathbf{a}))\\&\quad + \mathbb{P}(\mathbf{T}^\nu_\infty \leq \mathbf{s}|M^\nu(\mathbf{s}+\mathbf{a})^C)\mathbb{P}(M^\nu(\mathbf{s}+\mathbf{a})^C)\\&\to \mathbb{P}(\mathbf{Z} \leq \mathbf{s}) \qquad \text{as } \nu \to \infty,\end{aligned}$$

and the lemma is proved. $\quad\square$

By applying standard branching process theory (see, e.g., [15], if $R \leq 1$ (cf. Section 3), the total progeny of the Galton–Watson branching process is almost surely finite. In this case Lemma 4.1 is sufficient for studying $\mathbf{T}^\nu_\infty$. Alternatively, if $R > 1$, there is a positive probability $1 - \prod_{i=1}^{m}q_i^{a_i}$ that the Galton–Watson process produces an infinite number of progeny, where for $\mathbf{s} \in [\mathbf{0, 1}]$, $h_k(\mathbf{s}) = \mathbb{E}[\prod_{j=1}^m \exp((s_j - 1)\pi_j U_{k,j})]$ $(1 \le k \le m)$, and $\mathbf{q}$ is the unique solution in $[\mathbf{0, 1})$ of $\mathbf{q} = \mathbf{h}(\mathbf{q}) = (h_1(\mathbf{q}), h_2(\mathbf{q}), \ldots, h_m(\mathbf{q}))$.

Finally, using a lower-bound branching process approximation (see, e.g., [8, 9, 22]), it is straightforward to show that

$$\mathbb{P}\left(\lim_{\nu \to \infty} \bar{\mathbf{T}}_\infty^\nu = \mathbf{0}\right) = \prod_{i=1}^m q_i^{a_i} = 1 - \mathbb{P}\left(\lim_{\nu \to \infty} \bar{\mathbf{T}}_\infty^\nu = \boldsymbol{\tau}\right).$$

**5. Central limit theorem for the final size distribution.** In Section 3 we established a weak law of large numbers result for the final size distribution conditional upon there being a major outbreak ($\bar{\mathbf{T}}_\infty^\nu \xrightarrow{p} \boldsymbol{\tau}$ as $\nu \to \infty$). We shall extend the results of Section 3 to obtain a central limit theorem for $(\bar{\mathbf{T}}_\infty^\nu - \boldsymbol{\tau})\sqrt{N^\nu \Pi}$ as $\nu \to \infty$. We shall require the further constraints that for $1 \le i, j, k \le m$, $\sqrt{N^\nu}(N^\nu \mu_{ij}^\nu - \mu_{ij}) \to 0$, and there exists $0 \le \lambda_{ijk} < \infty$ such that $(N^\nu)^2 \lambda_{ijk}^\nu \to \lambda_{ijk}$ as $\nu \to \infty$. All standard epidemic models satisfy the above criteria (cf. Section 7). We also assume that $\sqrt{N^\nu}(\boldsymbol{\pi}^\nu - \boldsymbol{\pi}) \to \mathbf{0}$ as $\nu \to \infty$, with relaxations of this assumption considered in Section 6.

For $1 \le i \le m$ and $\mathbf{t} \ge \mathbf{0}$, let

$$Y_i^\nu(\mathbf{t}) = \sqrt{\tilde{N}_i^\nu}\left\{\frac{1}{N^\nu} X_i^\nu(\mathbf{t}) - \frac{\pi_i^\nu}{\pi_i} r_i^\nu(\mathbf{t})\right\}$$

and

$$\mathbf{Y}^\nu(\mathbf{t}) = \{\mathbf{X}^\nu(\mathbf{t})(N^\nu \Pi)^{-1} - \mathbf{r}^\nu(\mathbf{t})P^\nu\}\sqrt{N^\nu \Pi},$$

where $P^\nu = \Pi^\nu \Pi^{-1}$ and $\Pi^\nu = \mathrm{diag}(\pi_1^\nu, \pi_2^\nu, \ldots, \pi_m^\nu)$.

The general procedure in the literature for showing that $(\bar{\mathbf{T}}_\infty^\nu - \boldsymbol{\tau})\sqrt{N^\nu \Pi}$ converges in distribution to a multivariate Gaussian distribution is as follows. First, show that the finite-dimensional distributions of $\mathbf{Y}^\nu$ converge to the finite-dimensional distributions of a Gaussian process $\mathbf{Y}$. The next step is to show that $\{\mathbf{Y}^\nu\}$ are tight, and hence, $\mathbf{Y}^\nu \Rightarrow \mathbf{Y}$ as $\nu \to \infty$. Then by applying [11], Theorem 4.4, and pages 144–145, it follows that $\mathbf{Y}^\nu(\bar{\mathbf{T}}_\infty^\nu + \boldsymbol{\zeta}^\nu) \xrightarrow{D} \mathbf{Y}(\boldsymbol{\tau} + \boldsymbol{\zeta})$ as $\nu \to \infty$. Finally, it is shown that there exists an invertible $m \times m$ matrix $U$ such that $(\bar{\mathbf{T}}_\infty^\nu - \boldsymbol{\tau})\sqrt{N^\nu \Pi}$ and $\mathbf{Y}^\nu(\bar{\mathbf{T}}_\infty^\nu + \boldsymbol{\zeta}^\nu)U^{-1}$ have the same limiting distribution, namely, $\mathbf{Y}(\boldsymbol{\tau} + \boldsymbol{\zeta})U^{-1}$.

However, we shall give a simpler and more direct proof that $\mathbf{Y}^\nu(\bar{\mathbf{T}}_\infty^\nu + \boldsymbol{\zeta}^\nu) \xrightarrow{D} \mathbf{Y}(\boldsymbol{\tau} + \boldsymbol{\zeta})$ as $\nu \to \infty$. By [11], Theorem 4.1, it is sufficient to show that $\mathbf{Y}^\nu(\boldsymbol{\tau} + \boldsymbol{\zeta}) \xrightarrow{D} \mathbf{Y}(\boldsymbol{\tau} + \boldsymbol{\zeta})$ and $|\mathbf{Y}^\nu(\bar{\mathbf{T}}_\infty^\nu + \boldsymbol{\zeta}^\nu) - \mathbf{Y}^\nu(\boldsymbol{\tau} + \boldsymbol{\zeta})| \xrightarrow{p} 0$ as $\nu \to \infty$. These results are proved in Lemmas 5.3 and 5.4, respectively.

In proving a central limit theorem for $\{\mathbf{Y}^\nu(\mathbf{t})\}$, we shall prove that $\mathbf{Y}^\nu(\mathbf{t})$ has the same limiting distribution as $\hat{\mathbf{Y}}^\nu(\mathbf{t})P^\nu + \tilde{\mathbf{Y}}(\mathbf{t})$, where $\hat{\mathbf{Y}}^\nu(\mathbf{t})$ and $\tilde{\mathbf{Y}}(\mathbf{t})$ are independent and defined as follows. For $1 \le k \le m$, and $\mathbf{t} \ge 0$, let $Z_k^\nu(\mathbf{t}) =$ $\sum_{i=1}^{m} \sum_{j=1}^{[t_i \tilde{N}_i^\nu]} V_{(i,j),k}^\nu$, $z_k^\nu(\mathbf{t}) = \sum_{i=1}^{m} [t_i N^\nu \pi_i] \mu_{ik}^\nu$ and $z_k(\mathbf{t}) = \sum_{i=1}^{m} t_i \pi_i \mu_{ik}$, with $\mathbf{Z}^\nu(\mathbf{t}) = (Z_1^\nu(\mathbf{t}), Z_2^\nu(\mathbf{t}), \ldots, Z_m^\nu(\mathbf{t}))$, and $\mathbf{z}^\nu(\mathbf{t})$ and $\mathbf{z}(\mathbf{t})$ defined in the obvious fashion. For $\mathbf{b} \in \mathbb{R}^m$, let $\mathbf{g}(\mathbf{b}) = (\exp(-b_1), \exp(-b_2), \ldots, \exp(-b_m))$ and $G(\mathbf{b}) = \mathrm{diag}(\mathbf{g}(\mathbf{b}))$. Then

$$(5.1) \qquad \hat{\mathbf{Y}}^\nu(\mathbf{t}) = -\{\mathbf{g}(\mathbf{Z}^\nu(\mathbf{t})) - \mathbf{g}(\mathbf{z}^\nu(\mathbf{t}))\}\sqrt{N^\nu \Pi},$$

and $\hat{\mathbf{Y}}(\mathbf{t})$ is an $m$-dimensional Gaussian distribution with mean $\mathbf{0}$ and covariance matrix $(I - G(\mathbf{z}(\mathbf{t})))G(\mathbf{z}(\mathbf{t}))$. We shall therefore consider the limiting distribution of $\hat{\mathbf{Y}}^\nu(\mathbf{t})$ in Lemma 5.2, before considering $\mathbf{Y}^\nu$ in Lemma 5.3.

For $\mathbf{t} \geq \mathbf{0}$ and $1 \leq k \leq m$, let

$$\mathbf{W}_{\mathbf{t},k}^\nu = \sqrt{N^\nu} \sum_{j=1}^{[t_k \tilde{N}_k^\nu]} (\mathbf{V}_{(k,j)}^\nu - \boldsymbol{\mu}_k^\nu),$$

where $\boldsymbol{\mu}_k^\nu = (\mu_{k1}^\nu, \mu_{k2}^\nu, \ldots, \mu_{km}^\nu)$. Therefore $\mathbf{W}_{\mathbf{t},k}^\nu$ describes the fluctuations about the mean infectivity of the first $[t_k \tilde{N}_k^\nu]$ infectives of type $k$.

LEMMA 5.1. *Let* $\hat{\mathbf{W}}_{\mathbf{t}}$ *be an $m$-dimensional Gaussian distribution with mean $\mathbf{0}$ and covariance matrix* $\sqrt{\Pi}\{\sum_{k=1}^{m} t_k \pi_k \Lambda_k\}\sqrt{\Pi}$, *with* $\Lambda_k = (\lambda_{kij})$ $(1 \leq k \leq m)$. *Then for* $\mathbf{t} \geq \mathbf{0}$,

$$(5.2) \quad \hat{\mathbf{W}}_{\mathbf{t}}^\nu = \sum_{k=1}^{m} \{\mathbf{W}_{\mathbf{t},k}^\nu \sqrt{\Pi}\} \xrightarrow{D} \sum_{k=1}^{m} \{\mathbf{W}_{\mathbf{t},k} \sqrt{\Pi}\} = \hat{\mathbf{W}}_{\mathbf{t}}, \qquad \textit{say, as } \nu \to \infty.$$

PROOF.    The proof is trivial, and hence the details are omitted.    □

LEMMA 5.2. *For any* $\mathbf{t} \geq \mathbf{0}$,

$$\hat{\mathbf{Y}}^\nu(\mathbf{t}) \xrightarrow{D} \hat{\mathbf{W}}_{\mathbf{t}} G(\mathbf{z}(\mathbf{t})) \qquad \textit{as } \nu \to \infty.$$

PROOF.    By the mean value theorem, there exists $\mathbf{s}^\nu = (s_1^\nu, s_2^\nu, \ldots, s_m^\nu)$ lying between $\mathbf{Z}^\nu(\mathbf{t})$ and $\mathbf{z}^\nu(\mathbf{t})$ such that

$$\hat{\mathbf{Y}}^\nu(\mathbf{t}) = -\sum_{j=1}^{m} \left\{ (Z_j^\nu(\mathbf{t}) - z_j^\nu(\mathbf{t})) \frac{d}{dy_j} \mathbf{g}(\mathbf{y}) \Big|_{\mathbf{y}=\mathbf{s}^\nu} \right\} \sqrt{N^\nu \Pi^\nu}.$$

For $1 \leq j \leq m$,

$$\frac{d}{dy_j} \mathbf{g}(\mathbf{y}) \Big|_{\mathbf{y}=\mathbf{s}^\nu} = -(\delta_{1j} \exp(-s_1^\nu), \delta_{2j} \exp(-s_2^\nu), \ldots, \delta_{mj} \exp(-s_m^\nu)),$$

where $\delta_{ij}$ denotes the Kronecker delta. Thus

$$\hat{\mathbf{Y}}^\nu(\mathbf{t}) = (\mathbf{Z}^\nu(\mathbf{t}) - \mathbf{z}^\nu(\mathbf{t}))G(\mathbf{s}^\nu)\sqrt{N^\nu \Pi}$$
$$= \hat{\mathbf{W}}_{\mathbf{t}}^\nu G(\mathbf{s}^\nu).$$

By Lemma 5.1, $\mathbf{Z}^\nu(\mathbf{t}) \xrightarrow{p} \mathbf{z}(\mathbf{t})$, and hence, $\mathbf{s}^\nu \xrightarrow{p} \mathbf{z}(\mathbf{t})$ as $\nu \to \infty$. Therefore, by [11], Theorem 5.1, Corollary 2, $G(\mathbf{s}^\nu) \xrightarrow{p} G(\mathbf{z}(\mathbf{t}))$ as $\nu \to \infty$. The lemma then follows from Lemma 5.1. $\square$

LEMMA 5.3. *For any* $\mathbf{t} \geq \mathbf{0}$,

$$\mathbf{Y}^\nu(\mathbf{t}) \xrightarrow{D} \mathbf{Y}(\mathbf{t}) \qquad as \ \nu \to \infty,$$

*where* $\mathbf{Y}(\mathbf{t})$ *is an m-dimensional Gaussian distribution with mean* $\mathbf{0}$ *and covariance matrix*

$$\tilde{\Xi}(\mathbf{t}) = \{I - G(\mathbf{z}(\mathbf{t}))\}G(\mathbf{z}(\mathbf{t})) + G(\mathbf{z}(\mathbf{t}))\sqrt{\Pi}\left\{\sum_{k=1}^{m} t_k \pi_k \Lambda_k\right\}\sqrt{\Pi}G(\mathbf{z}(\mathbf{t})).$$

PROOF. Since $P^\nu \to I$ as $\nu \to \infty$, $\mathbf{Y}(\mathbf{t})$ is the limiting distribution of $\hat{\mathbf{Y}}^\nu(\mathbf{t})P^\nu + \tilde{\mathbf{Y}}(\mathbf{t})$ as $\nu \to \infty$.

For $\boldsymbol{\gamma} \in \mathbb{R}^m$, let $Z^\nu_{\boldsymbol{\gamma}}(\mathbf{t}) = \sum_{l=1}^{m} \gamma_l Y^\nu_l(\mathbf{t})$, and for $s \in \mathbb{R}$, let

$$f_\nu(s; \mathbf{t}; \boldsymbol{\gamma}) = \mathbb{E}[\exp(isZ^\nu_{\boldsymbol{\gamma}}(\mathbf{t}))].$$

Let

$$\tilde{f}_\nu(s; \mathbf{t}; \boldsymbol{\gamma}) = \prod_{k=1}^{m} \exp\left(-\frac{s^2 \gamma_k^2}{2} r_k(\mathbf{t})\{1 - r_k(\mathbf{t})\}\right)$$

$$\times \mathbb{E}_{\tilde{\mathbf{V}}^\nu}\left[\prod_{k=1}^{m} \exp\left(is\gamma_k \frac{N^\nu_k}{\sqrt{\tilde{N}^\nu_k}}\{r^\nu_k(\mathbf{t}|\tilde{\mathbf{V}}^\nu) - r^\nu_k(\mathbf{t})\}\right)\right]$$

$$= \sum_{l=1}^{m} \gamma_l \left(\tilde{Y}_l(\mathbf{t}) + \frac{\pi^\nu_l}{\pi_l}\hat{Y}^\nu_l(\mathbf{t})\right).$$

Therefore, using the Cramér–Wold device and characteristic functions, the lemma follows by proving that

(5.3) $\qquad |f_\nu(s; \mathbf{t}; \boldsymbol{\gamma}) - \tilde{f}_\nu(s; \mathbf{t}; \boldsymbol{\gamma})| \to 0 \qquad$ as $\nu \to \infty$.

Conditional upon $\tilde{\mathbf{V}}^\nu$, $\chi^\nu_{i,j}(\mathbf{t})$ and $\chi^\nu_{i',j'}(\mathbf{t})$ are independent, and so,

$$f_\nu(s; \mathbf{t}; \boldsymbol{\gamma}) = \mathbb{E}_{\tilde{\mathbf{V}}^\nu}[\mathbb{E}_{Z^\nu_{\boldsymbol{\gamma}}}[\exp(isZ^\nu_{\boldsymbol{\gamma}}(\mathbf{t}))|\tilde{\mathbf{V}}^\nu]]$$

$$= \mathbb{E}_{\tilde{\mathbf{V}}^\nu}\left[\prod_{k=1}^{m}\prod_{j=1}^{N^\nu_k} \mathbb{E}_{\chi^\nu_{k,j}}\left[\exp\left(is\gamma_k \frac{1}{\sqrt{\tilde{N}^\nu_k}}(\chi^\nu_{k,j}(\mathbf{t}) - r^\nu_k(\mathbf{t}))\right)\Big|\tilde{\mathbf{V}}^\nu\right]\right].$$

For $1 \leq k \leq m$, and $1 \leq j \leq N^\nu_k$,

$$\mathbb{E}_{\chi^\nu_{k,j}}\left[\exp\left(is\gamma_k \frac{1}{\sqrt{\tilde{N}^\nu_k}}(\chi^\nu_{k,j}(\mathbf{t}) - r^\nu_k(\mathbf{t}))\right)\Big|\tilde{\mathbf{V}}^\nu\right]$$

$$= \exp\left(is\frac{\gamma_k}{\sqrt{\tilde{N}_k^\nu}}(r_k^\nu(\mathbf{t}|\tilde{\mathbf{V}}^\nu) - r_k^\nu(\mathbf{t}))\right)$$
$$\times \mathbb{E}_{\chi_{k,j}^\nu}\left[\exp\left(is\frac{\gamma_k}{\sqrt{\tilde{N}_k^\nu}}\{\chi_{k,j}^\nu(\mathbf{t}) - r_k^\nu(\mathbf{t}|\tilde{\mathbf{V}}^\nu)\}\right)\right]\Big|\tilde{\mathbf{V}}^\nu\right].$$

Let
$$\tilde{f}_\nu^1(s;\mathbf{t};\boldsymbol{\gamma}) = \mathbb{E}_{\tilde{\mathbf{V}}^\nu}\left[\prod_{k=1}^m\left\{\exp\left(is\frac{N_k^\nu}{\sqrt{\tilde{N}_k^\nu}}\gamma_k\{r_k^\nu(\mathbf{t}|\tilde{\mathbf{V}}^\nu) - r_k^\nu(\mathbf{t})\}\right)\right.\right.$$
$$\left.\left.\times\left(1 - \frac{s^2\gamma_k^2}{2\tilde{N}_k^\nu}r_k^\nu(\mathbf{t}|\tilde{\mathbf{V}}^\nu)\{1 - r_k^\nu(\mathbf{t}|\tilde{\mathbf{V}}^\nu)\}\right)^{N_k^\nu}\right\}\right]$$

where $\text{var}(\chi_{k,j}|\tilde{\mathbf{V}}^\nu) = r_k^\nu(\mathbf{t}|\tilde{\mathbf{V}}^\nu)\{1 - r_k^\nu(\mathbf{t}|\tilde{\mathbf{V}}^\nu)\}$ $(1 \le k \le m, 1 \le j \le N_k^\nu)$.

For any $s \in \mathbb{R}$, and for any random variable $W$ such that $\mathbb{E}[|W|^3] < \infty$, there exists $n_0 \in \mathbb{N}$ such that for all $n \ge n_0$,

$$\begin{aligned}(5.4)\quad &\left|\mathbb{E}\left[\exp\left(is\frac{1}{\sqrt{n}}W\right)\right] - \left(1 + \frac{is}{\sqrt{n}}\mathbb{E}[W] - \frac{s^2}{2n}\mathbb{E}[W^2]\right)\right| \\ &\le \left(\frac{|s|}{\sqrt{n}}\right)^3\mathbb{E}[|W|^3].\end{aligned}$$

Therefore for all $1 \le k \le m$, $1 \le j \le N_k^\nu$ and $\mathbf{t} \ge \mathbf{0}$, it follows from (5.4) that

$$|f_\nu(s;\mathbf{t};\boldsymbol{\gamma}) - \tilde{f}_\nu^1(s;\mathbf{t};\boldsymbol{\gamma})| \le \sum_{k=1}^m N_k^\nu\left(\frac{|s\gamma_k|}{\sqrt{\tilde{N}_k^\nu}}\right)^3 \to 0 \qquad \text{as } \nu \to \infty.$$

Let
$$\tilde{f}_\nu^2(s;\mathbf{t};\boldsymbol{\gamma}) = \mathbb{E}_{\tilde{\mathbf{V}}^\nu}\left[\prod_{k=1}^m\left\{\exp\left(is\frac{N_k^\nu}{\sqrt{\tilde{N}_k^\nu}}\gamma_k\{r_k^\nu(\mathbf{t}|\tilde{\mathbf{V}}^\nu) - r_k^\nu(\mathbf{t})\}\right)\right.\right.$$
$$\left.\left.\times\exp\left(-\frac{s^2}{2}\frac{N_k^\nu}{\tilde{N}_k^\nu}\gamma_k^2 r_k^\nu(\mathbf{t}|\tilde{\mathbf{V}}^\nu)\{1 - r_k^\nu(\mathbf{t}|\tilde{\mathbf{V}}^\nu)\}\right)\right)\right\}\right].$$

It is trivial, using [12], page 94, Lemma 4.3, to show that

$$|\tilde{f}_\nu^1(s;\mathbf{t};\boldsymbol{\gamma}) - \tilde{f}_\nu^2(s;\mathbf{t};\boldsymbol{\gamma})| \le \sum_{k=1}^m N_k^\nu\left(\frac{s^2\gamma_k^2}{2\tilde{N}_k^\nu}\right)^2 \to 0 \qquad \text{as } \nu \to \infty.$$

For $1 \le k \le m$,

$$(5.5)\quad \left|r_k^\nu(\mathbf{t}|\tilde{\mathbf{V}}^\nu) - \left(1 - \exp\left(-\sum_{i=1}^m\sum_{j=1}^{[t_iN_i^\nu]}V_{(i,j),k}^\nu\right)\right)\right| \le \sum_{i=1}^m\sum_{j=1}^{[t_iN_i^\nu]}(V_{(i,j),k}^\nu)^2,$$

with $\mathbb{E}[\sum_{i=1}^{m}\sum_{j=1}^{[t_i N_i^{\nu}]}(V_{(i,j),k}^{\nu})^2] = \sum_{i=1}^{m}[t_i N_i^{\nu}](\lambda_{ikk}^{\nu} + (\mu_{ik}^{\nu})^2) \to 0$ as $\nu \to \infty$. Therefore, by Lemma 5.1, and [11], Theorem 5.1, Corollary 2,

$$1 - \exp\left(-\sum_{i=1}^{m}\sum_{j=1}^{[t_i N_i^{\nu}]}V_{(i,j),k}^{\nu}\right) \overset{p}{\longrightarrow} 1 - \exp\left(-\sum_{i=1}^{m}t_i\pi_i\mu_{ik}\right) = r_k(\mathbf{t}),$$

and so, $r_k^{\nu}(\mathbf{t}|\tilde{\mathbf{V}}^{\nu}) \overset{p}{\longrightarrow} r_k(\mathbf{t})$ as $\nu \to \infty$. For all $1 \le k \le m$, $N_k^{\nu}/\tilde{N}_k^{\nu} \to 1$ as $\nu \to \infty$ and for all $\nu \ge 1$, $|\tilde{f}_{\nu}^2(s;\mathbf{t};\boldsymbol{\gamma})|, |\tilde{f}_{\nu}(s;\mathbf{t};\boldsymbol{\gamma})| \le 1$. Therefore

$$|\tilde{f}_{\nu}^2(s;\mathbf{t};\boldsymbol{\gamma}) - \tilde{f}_{\nu}(s;\mathbf{t};\boldsymbol{\gamma})| \to 0 \qquad \text{as } \nu \to \infty.$$

Thus (5.3) is proved. □

For $1 \le k \le m$, $\exp(-z_k(\boldsymbol{\tau} + \boldsymbol{\zeta})) = \sigma_k$, and so,

$$(5.6) \qquad G(\mathbf{z}(\boldsymbol{\tau} + \boldsymbol{\zeta})) = \Sigma \equiv \operatorname{diag}(\sigma_1, \sigma_2, \dots, \sigma_m).$$

Therefore, $\mathbf{Y}(\boldsymbol{\tau} + \boldsymbol{\zeta}) \sim N(\mathbf{0}, \Xi)$, where

$$(5.7) \qquad \Xi = \{I - \Sigma\}\Sigma + \Sigma\sqrt{\Pi}\left\{\sum_{k=1}^{m}(\tau_k + \zeta_k)\pi_k\Lambda_k\right\}\sqrt{\Pi}\Sigma.$$

The next step is to show that $\mathbf{Y}^{\nu}(\bar{\mathbf{T}}_{\infty}^{\nu} + \boldsymbol{\zeta}^{\nu})$ and $\mathbf{Y}^{\nu}(\boldsymbol{\tau} + \boldsymbol{\zeta})$ have the same limiting distribution.

LEMMA 5.4. *Suppose that $\bar{\mathbf{T}}_{\infty}^{\nu} \overset{p}{\longrightarrow} \boldsymbol{\tau}$. Then*

$$\mathbf{Y}^{\nu}(\bar{\mathbf{T}}_{\infty}^{\nu} + \boldsymbol{\zeta}^{\nu}) \overset{D}{\longrightarrow} \mathbf{Y}(\boldsymbol{\tau} + \boldsymbol{\zeta}) \qquad \text{as } \nu \to \infty.$$

PROOF. By Lemma 5.3, $\mathbf{Y}^{\nu}(\boldsymbol{\tau} + \boldsymbol{\zeta}) \overset{D}{\longrightarrow} \mathbf{Y}(\boldsymbol{\tau} + \boldsymbol{\zeta})$ as $\nu \to \infty$. Therefore, by [11], Theorem 4.1, to prove the lemma it is sufficient to show that

$$|\mathbf{Y}^{\nu}(\bar{\mathbf{T}}_{\infty}^{\nu} + \boldsymbol{\zeta}^{\nu}) - \mathbf{Y}^{\nu}(\boldsymbol{\tau} + \boldsymbol{\zeta})| \overset{p}{\longrightarrow} 0 \qquad \text{as } \nu \to \infty.$$

Fix $\varepsilon > 0$. Then for $1 \le i \le m$, by Chebyshev's inequality

$$\mathbb{P}(|Y_i^{\nu}(\bar{\mathbf{T}}_{\infty}^{\nu} + \boldsymbol{\zeta}^{\nu}) - Y_i^{\nu}(\boldsymbol{\tau} + \boldsymbol{\zeta})| > \varepsilon)$$
$$= \sum_{\mathbf{k}^{\nu}}\mathbb{P}(|Y_i^{\nu}(\bar{\mathbf{T}}_{\infty}^{\nu} + \boldsymbol{\zeta}^{\nu}) - Y_i^{\nu}(\boldsymbol{\tau} + \boldsymbol{\zeta})| > \varepsilon|\bar{\mathbf{T}}_{\infty}^{\nu} + \boldsymbol{\zeta}^{\nu} = \mathbf{k}^{\nu})\mathbb{P}(\bar{\mathbf{T}}_{\infty}^{\nu} + \boldsymbol{\zeta}^{\nu} = \mathbf{k}^{\nu})$$
$$\le \sum_{\mathbf{k}^{\nu}}\frac{1}{\varepsilon^2}\operatorname{var}(Y_i^{\nu}(\bar{\mathbf{T}}_{\infty}^{\nu} + \boldsymbol{\zeta}^{\nu}) - Y_i^{\nu}(\boldsymbol{\tau} + \boldsymbol{\zeta})|\bar{\mathbf{T}}_{\infty}^{\nu} + \boldsymbol{\zeta}^{\nu} = \mathbf{k}^{\nu})\mathbb{P}(\bar{\mathbf{T}}_{\infty}^{\nu} + \boldsymbol{\zeta}^{\nu} = \mathbf{k}^{\nu})$$
$$= \frac{1}{\varepsilon^2}\mathbb{E}[\operatorname{var}(Y_i^{\nu}(\bar{\mathbf{T}}_{\infty}^{\nu} + \boldsymbol{\zeta}^{\nu}) - Y_i^{\nu}(\boldsymbol{\tau} + \boldsymbol{\zeta})|\bar{\mathbf{T}}_{\infty}^{\nu} + \boldsymbol{\zeta}^{\nu})].$$

Thus

$$\mathbb{P}(|Y_i^\nu(\bar{\mathbf{T}}_\infty^\nu + \boldsymbol{\zeta}^\nu) - Y_i^\nu(\boldsymbol{\tau} + \boldsymbol{\zeta})| > \varepsilon)$$

$$\leq \frac{1}{\varepsilon^2 \bar{N}_i^\nu} \sum_{j=1}^{N_i^\nu} \sum_{l=1}^{N_i^\nu} \mathbb{E}[\mathrm{cov}(\chi_{i,j}^\nu(\bar{\mathbf{T}}_\infty^\nu + \boldsymbol{\zeta}^\nu) - \chi_{i,j}^\nu(\boldsymbol{\tau} + \boldsymbol{\zeta}),$$

$$\chi_{i,l}^\nu(\bar{\mathbf{T}}_\infty^\nu + \boldsymbol{\zeta}^\nu) - \chi_{i,l}^\nu(\boldsymbol{\tau} + \boldsymbol{\zeta})|\bar{\mathbf{T}}_\infty^\nu + \boldsymbol{\zeta}^\nu)]$$

$$(5.8) \qquad = \frac{\pi_i^\nu}{\pi_i \varepsilon^2} \{\mathbb{E}[\mathrm{var}(\chi_{i,1}^\nu(\bar{\mathbf{T}}_\infty^\nu + \boldsymbol{\zeta}^\nu) - \chi_{i,1}^\nu(\boldsymbol{\tau} + \boldsymbol{\zeta})|\bar{\mathbf{T}}_\infty^\nu + \boldsymbol{\zeta}^\nu)]$$

$$+ (N_i^\nu - 1)\mathbb{E}[\mathrm{cov}(\chi_{i,1}^\nu(\bar{\mathbf{T}}_\infty^\nu + \boldsymbol{\zeta}^\nu) - \chi_{i,1}^\nu(\boldsymbol{\tau} + \boldsymbol{\zeta}),$$

$$\chi_{i,2}^\nu(\bar{\mathbf{T}}_\infty^\nu + \boldsymbol{\zeta}^\nu) - \chi_{i,2}^\nu(\boldsymbol{\tau} + \boldsymbol{\zeta})|\bar{\mathbf{T}}_\infty^\nu + \boldsymbol{\zeta}^\nu)]\}.$$

The lemma follows by proving that the right-hand side of (5.8) converges to 0 as $\nu \to \infty$.

First,

$$(5.9) \qquad \mathbb{E}[\mathrm{var}(\chi_{i,1}^\nu(\bar{\mathbf{T}}_\infty^\nu + \boldsymbol{\zeta}^\nu) - \chi_{i,1}^\nu(\boldsymbol{\tau} + \boldsymbol{\zeta})|\bar{\mathbf{T}}_\infty^\nu + \boldsymbol{\zeta}^\nu)]$$

$$\leq \mathbb{E}[|r_i^\nu(\bar{\mathbf{T}}_\infty^\nu + \boldsymbol{\zeta}^\nu) - r_i^\nu(\boldsymbol{\tau} + \boldsymbol{\zeta})|].$$

Since $\bar{\mathbf{T}}_\infty^\nu \xrightarrow{p} \boldsymbol{\tau}$, it follows that $|r_i^\nu(\bar{\mathbf{T}}_\infty^\nu + \boldsymbol{\zeta}^\nu) - r_i^\nu(\boldsymbol{\tau} + \boldsymbol{\zeta})| \xrightarrow{p} 0$ as $\nu \to \infty$. Therefore, since $r_i^\nu(\cdot)$ is bounded, the right-hand side of (5.9) converges to 0 as $\nu \to \infty$.

By (3.3) of Lemma 3.1,

$$N_i^\nu |\mathrm{cov}(\chi_{i,1}^\nu(\bar{\mathbf{T}}_\infty^\nu + \boldsymbol{\zeta}^\nu) - \chi_{i,1}^\nu(\boldsymbol{\tau} + \boldsymbol{\zeta}),$$

$$\chi_{i,2}^\nu(\bar{\mathbf{T}}_\infty^\nu + \boldsymbol{\zeta}^\nu) - \chi_{i,2}^\nu(\boldsymbol{\tau} + \boldsymbol{\zeta})|\bar{\mathbf{T}}_\infty^\nu + \boldsymbol{\zeta}^\nu)|$$

$$(5.10) \qquad \leq 6N_i^\nu \left\{ \sum_{k=1}^m ([(\bar{T}_{\infty,k}^\nu \vee \tau_k)\tilde{N}_k^\nu] - [(\bar{T}_{\infty,k}^\nu \wedge \tau_k)\tilde{N}_k^\nu])\mu_{ki}^\nu \right\}$$

$$\times \left\{ \sum_{k=1}^m (\bar{T}_{\infty,k}^\nu \wedge \tau_k)\tilde{N}_k^\nu \lambda_{kii}^\nu \right\}$$

$$+ 2N_i^\nu \sum_{k=1}^m ([(\bar{T}_{\infty,k}^\nu \vee \tau_k)\tilde{N}_k^\nu] - [(\bar{T}_{\infty,k}^\nu \wedge \tau_k)\tilde{N}_k^\nu])\lambda_{kii}^\nu.$$

It is trivial to show that for all $1 \leq k \leq m$,

$$\frac{1}{N^\nu}\{[(\bar{T}_{\infty,k}^\nu \vee \tau_k)\tilde{N}_k^\nu] - [(\bar{T}_{\infty,k}^\nu \wedge \tau_k)\tilde{N}_k^\nu]\} \xrightarrow{p} 0 \qquad \text{as } \nu \to \infty.$$

Therefore, since $N^\nu \mu_{ki}^\nu \to \mu_{ki}$ and $(N^\nu)^2 \lambda_{kii}^\nu \to \lambda_{kii}$ as $\nu \to \infty$, the right-hand side of (5.10) converges in probability to 0 as $\nu \to \infty$. Finally, since

$\mathbf{0} \le \bar{\mathbf{T}}_\infty^\nu \le \mathbf{1}$,

$$N_i^\nu \mathbb{E}[|\operatorname{cov}(\chi_{i,1}^\nu(\bar{\mathbf{T}}_\infty^\nu + \boldsymbol{\zeta}^\nu) - \chi_{i,1}(\boldsymbol{\tau} + \boldsymbol{\zeta}),$$
$$\chi_{i,2}^\nu(\bar{\mathbf{T}}_\infty^\nu + \boldsymbol{\zeta}^\nu) - \chi_{i,2}(\boldsymbol{\tau} + \boldsymbol{\zeta})|\bar{\mathbf{T}}_\infty^\nu + \boldsymbol{\zeta}^\nu)|] \to 0$$

as $\nu \to \infty$.  $\square$

We are now in position to derive a Gaussian limit for $(\bar{\mathbf{T}}_\infty^\nu - \boldsymbol{\tau})\sqrt{N^\nu \Pi}$.

THEOREM 5.5. *Let $U$ be the $m \times m$ matrix with entries $u_{ij} = \delta_{ij} - \sqrt{\pi_i \pi_j} \mu_{ij} \sigma_j$. Then $U$ is invertible. Furthermore, if $(\boldsymbol{\zeta}^\nu - \boldsymbol{\zeta})\sqrt{N^\nu \Pi^\nu} \to \mathbf{0}$ and $\bar{\mathbf{T}}_\infty^\nu \xrightarrow{p} \boldsymbol{\tau}$ as $\nu \to \infty$, then*

$$(\bar{\mathbf{T}}_\infty^\nu - \boldsymbol{\tau})\sqrt{N^\nu \Pi} \xrightarrow{D} N(\mathbf{0}, (U^T)^{-1}\Xi U^{-1}) \qquad \text{as } \nu \to \infty,$$

*where $\Xi$ is defined in* (5.7).

PROOF.  Note that

$$(\bar{\mathbf{T}}_\infty^\nu - \boldsymbol{\tau})\sqrt{N^\nu \Pi}$$
$$= (\mathbf{X}^\nu(\bar{\mathbf{T}}_\infty^\nu + \boldsymbol{\zeta}^\nu)(N^\nu \Pi)^{-1} - \mathbf{r}(\boldsymbol{\tau} + \boldsymbol{\zeta}))\sqrt{N^\nu \Pi}$$
$$= \mathbf{Y}^\nu(\bar{\mathbf{T}}_\infty^\nu + \boldsymbol{\zeta}^\nu) + \{\mathbf{r}^\nu(\bar{\mathbf{T}}_\infty^\nu + \boldsymbol{\zeta}^\nu)P^\nu - \mathbf{r}(\bar{\mathbf{T}}_\infty^\nu + \boldsymbol{\zeta}^\nu)\}\sqrt{N^\nu \Pi}$$
$$\quad + \{\mathbf{r}(\bar{\mathbf{T}}_\infty^\nu + \boldsymbol{\zeta}^\nu) - \mathbf{r}(\boldsymbol{\tau} + \boldsymbol{\zeta})\}\sqrt{N^\nu \Pi}.$$

Hence,

$$(\{\bar{\mathbf{T}}_\infty^\nu + \boldsymbol{\zeta}^\nu - \mathbf{r}(\bar{\mathbf{T}}_\infty^\nu + \boldsymbol{\zeta}^\nu)\} - \{\boldsymbol{\tau} + \boldsymbol{\zeta} - \mathbf{r}(\boldsymbol{\tau} + \boldsymbol{\zeta})\})\sqrt{N^\nu \Pi}$$
$$= \mathbf{Y}^\nu(\bar{\mathbf{T}}_\infty^\nu + \boldsymbol{\zeta}^\nu) + \{\mathbf{r}^\nu(\bar{\mathbf{T}}_\infty^\nu + \boldsymbol{\zeta}^\nu)P^\nu - \mathbf{r}(\bar{\mathbf{T}}_\infty^\nu + \boldsymbol{\zeta}^\nu)\}\sqrt{N^\nu \Pi}$$
$$\quad + \{\boldsymbol{\zeta}^\nu - \boldsymbol{\zeta}\}\sqrt{N^\nu \Pi}.$$

By the mean value theorem,

$$(\{\bar{\mathbf{T}}_\infty^\nu + \boldsymbol{\zeta}^\nu - \mathbf{r}(\bar{\mathbf{T}}_\infty^\nu + \boldsymbol{\zeta}^\nu)\} - \{\boldsymbol{\tau} + \boldsymbol{\zeta} - \mathbf{r}(\boldsymbol{\tau} + \boldsymbol{\zeta})\})\sqrt{N^\nu \Pi}$$
$$= \left(\sum_{j=1}^m ((\bar{T}_{\infty,j}^\nu + \zeta_j^\nu) - (\tau_j + \zeta_j))\frac{d}{dt_j}(\mathbf{t} - \mathbf{r}(\mathbf{t}))\Big|_{\mathbf{t}=\mathbf{s}^\nu}\right)\sqrt{N^\nu \Pi},$$

where $\mathbf{s}^\nu$ lies between $\bar{\mathbf{T}}_\infty^\nu + \boldsymbol{\zeta}^\nu$ and $\boldsymbol{\tau} + \boldsymbol{\zeta}$.

For $\mathbf{t} \ge \mathbf{0}$, let $S(\mathbf{t})$ and $U(\mathbf{t})$ denote the $m \times m$ matrices with entries, $s_{ij}(\mathbf{t}) = \delta_{ij} - \pi_i \mu_{ij}(1 - r_j(\mathbf{t}))$ and $u_{ij}(\mathbf{t}) = \delta_{ij} - \sqrt{\pi_i \pi_j}\mu_{ij}(1 - r_j(\mathbf{t}))$. Therefore

$$(5.11) \quad (\bar{\mathbf{T}}_\infty^\nu - \boldsymbol{\tau})S(\mathbf{s}^\nu)\sqrt{N^\nu \Pi} = (\bar{\mathbf{T}}_\infty^\nu - \boldsymbol{\tau})\sqrt{N^\nu \Pi}U(\mathbf{s}^\nu)$$
$$= \mathbf{Y}^\nu(\bar{\mathbf{T}}_\infty^\nu + \boldsymbol{\zeta}^\nu) + \mathbf{A}_\nu, \qquad \text{say,}$$

where $\mathbf{A}_\nu \xrightarrow{p} 0$ as $\nu \to \infty$. Thus by [11], Theorem 4.1, $(\bar{\mathbf{T}}_\infty^\nu - \boldsymbol{\tau})\sqrt{N^\nu \Pi}$ and $\mathbf{Y}^\nu(\bar{\mathbf{T}}_\infty^\nu + \boldsymbol{\zeta})U(\mathbf{s}^\nu)^{-1}$ have the same limiting distribution as $\nu \to \infty$, provided that $\lim_{\nu \to \infty} U(\mathbf{s}^\nu)$ exists.

Since $U(\cdot)$ is continuous, and $\mathbf{s}^\nu \xrightarrow{p} \boldsymbol{\tau} + \boldsymbol{\zeta}$ as $\nu \to \infty$, it follows by [11], Theorem 5.1, Corollary 2, that

$$U(\mathbf{s}^\nu) \xrightarrow{p} U(\boldsymbol{\tau} + \boldsymbol{\zeta}) \equiv U \qquad \text{as } \nu \to \infty.$$

By [20], page 575, $U$ is invertible. Therefore, by combining Lemma 5.4 and [11], Theorem 4.4, we have that

$$(5.12) \quad \mathbf{Y}^\nu(\bar{\mathbf{T}}_\infty^\nu + \boldsymbol{\zeta}^\nu)U(\mathbf{s}^\nu)^{-1} \xrightarrow{D} \mathbf{Y}(\boldsymbol{\tau} + \boldsymbol{\zeta})U^{-1} \sim N(\mathbf{0}, (U^T)^{-1}\Xi U^{-1})$$
$$\text{as } \nu \to \infty,$$

and the theorem is proved. $\quad\square$

## 6. Random group allocation.

There is an interesting variant of the above model in which individuals are independently and randomly allocated to a group according to a probability distribution, $\boldsymbol{\pi}$. The main motivation for studying this variation of the model is that it enables us to consider a wider range of epidemics upon random graphs in Section 7.

Let $\mathbf{S}$ denote an $m$-dimensional Gaussian random variable with mean $\mathbf{0}$ and covariance matrix $\tilde{\Upsilon} = \Pi - \boldsymbol{\pi}^T\boldsymbol{\pi}$. Then $\sqrt{N^\nu}(\boldsymbol{\pi}^\nu - \boldsymbol{\pi}) \xrightarrow{D} \mathbf{S}$ as $\nu \to \infty$.

To assist in the analysis of the model, we define for $1 \le i \le m$, $\tilde{X}_i^\nu(\mathbf{t}) = \sum_{j=1}^{\tilde{N}_i^\nu} \chi_{i,j}^\nu(\mathbf{t})$ and $\tilde{Y}_i^\nu = \sqrt{\tilde{N}_i^\nu}(\frac{1}{\tilde{N}_i^\nu}\tilde{X}_i^\nu(\mathbf{t}) - r_i^\nu(\mathbf{t}))$, with $\tilde{N}_i^\nu = N^\nu \pi_i$ (cf. Section 3). Also $\tilde{\mathbf{X}}^\nu(\mathbf{t})$ and $\tilde{\mathbf{Y}}^\nu(\mathbf{t})$ are defined in the obvious fashion. Note that $\mathbf{N}^\nu$ is a random vector, while $\tilde{\mathbf{N}}^\nu$ is deterministic.

The results proved in Sections 3 and 4 are valid for the above model with only trivial modifications of the proofs required. However, Theorem 5.5 requires minor alterations for the current model. We begin by considering $\{\mathbf{Y}^\nu(\bar{\mathbf{T}}_\infty^\nu + \boldsymbol{\zeta}^\nu)\}$. We then prove Theorem 6.2 which is an adaption of Theorem 5.5.

Lemmas 5.1–5.4 hold with $\{\tilde{\mathbf{Y}}^\nu(\bar{\mathbf{T}}_\infty^\nu + \boldsymbol{\zeta}^\nu)\}$ in place of $\{\mathbf{Y}^\nu(\bar{\mathbf{T}}_\infty^\nu + \boldsymbol{\zeta}^\nu)\}$. Therefore we shall require Lemma 6.1 to utilize these results for the current model.

LEMMA 6.1.

$$|\mathbf{Y}^\nu(\bar{\mathbf{T}}_\infty^\nu + \boldsymbol{\zeta}^\nu) - \tilde{\mathbf{Y}}^\nu(\bar{\mathbf{T}}_\infty^\nu + \boldsymbol{\zeta}^\nu)| \xrightarrow{p} 0 \qquad as \ \nu \to \infty.$$

PROOF. Fix $1 \le i \le m$ and $\varepsilon > 0$. Then by applying Chebyshev's inequality in a similar manner to Lemma 5.4, we have that

$$\mathbb{P}(|\tilde{Y}_i^\nu(\bar{\mathbf{T}}_\infty^\nu + \boldsymbol{\zeta}^\nu) - Y_i^\nu(\bar{\mathbf{T}}_\infty^\nu + \boldsymbol{\zeta}^\nu)| > \varepsilon)$$

$$\leq \frac{1}{\varepsilon^2} \mathbb{E}[\operatorname{var}(\tilde{Y}_i^\nu(\bar{\mathbf{T}}_\infty^\nu + \boldsymbol{\zeta}^\nu) - Y_i^\nu(\bar{\mathbf{T}}_\infty^\nu + \boldsymbol{\zeta}^\nu)|\mathbf{N}^\nu)]$$

$$= \frac{1}{\varepsilon^2} \mathbb{E}\left[\frac{1}{\tilde{N}_i^\nu} \operatorname{var}\left(\sum_{j=1}^{\tilde{N}_i^\nu} \chi_{i,j}^\nu(\bar{\mathbf{T}}_\infty^\nu + \boldsymbol{\zeta}^\nu) - \sum_{j=1}^{N_i^\nu} \chi_{i,j}^\nu(\bar{\mathbf{T}}_\infty^\nu + \boldsymbol{\zeta}^\nu)\Big|\mathbf{N}^\nu\right)\right]$$

$$= \frac{1}{\varepsilon^2 \tilde{N}_i^\nu} \mathbb{E}\left[\sum_{j=(\tilde{N}_i^\nu \wedge N_i^\nu)+1}^{\tilde{N}_i^\nu \vee N_i^\nu} \sum_{l=(\tilde{N}_i^\nu \wedge N_i^\nu)+1}^{\tilde{N}_i^\nu \vee N_i^\nu} \operatorname{cov}(\chi_{i,j}^\nu(\bar{\mathbf{T}}_\infty^\nu + \boldsymbol{\zeta}^\nu), \chi_{i,l}^\nu(\bar{\mathbf{T}}_\infty^\nu + \boldsymbol{\zeta}^\nu))\right].$$

Since $\bar{\mathbf{T}}_\infty^\nu$ is bounded above by $\mathbf{1}$, it follows from (3.2) that

$$\begin{aligned}
(6.1) \qquad &\mathbb{P}(|\tilde{Y}_i^\nu(\bar{\mathbf{T}}_\infty^\nu + \boldsymbol{\zeta}^\nu) - Y_i^\nu(\bar{\mathbf{T}}_\infty^\nu + \boldsymbol{\zeta}^\nu)| > \varepsilon) \\
&\leq \frac{1}{\varepsilon^2 \tilde{N}_i^\nu} \mathbb{E}\left[|\tilde{N}_i^\nu - N_i^\nu|\left\{1 + |\tilde{N}_i^\nu - N_i^\nu| \sum_{k=1}^m \tilde{N}_k^\nu \lambda_{kii}^\nu\right\}\right] \\
&\leq \frac{1}{\varepsilon^2}\left\{\mathbb{E}\left[\frac{|\pi_i^\nu - \pi_i|}{\pi_i}\right] + \sup_{1 \leq k \leq m}(N^\nu)^2 \lambda_{kii}^\nu \mathbb{E}\left[\frac{(\pi_i^\nu - \pi_i)^2}{\pi_i}\right]\right\}.
\end{aligned}$$

Thus the right-hand side of (6.1) converges to 0 as $\nu \to \infty$, and the lemma is proved. $\square$

THEOREM 6.2.  *Suppose that* $(\boldsymbol{\zeta}^\nu - \boldsymbol{\zeta})\sqrt{N^\nu\Pi} \xrightarrow{p} \mathbf{0}$ *and* $\bar{\mathbf{T}}_\infty^\nu \xrightarrow{p} \boldsymbol{\tau}$ *as* $\nu \to \infty$. *Then*

$$(\bar{\mathbf{T}}_\infty^\nu - \boldsymbol{\tau})\sqrt{N^\nu\Pi} \xrightarrow{D} N(\mathbf{0}, (U^T)^{-1}\{\Xi + \Upsilon\}U^{-1}) \qquad as\ \nu \to \infty,$$

*where* $\Xi$ *and* $U$ *are defined in* (5.7) *and Theorem* 5.5, *respectively, and*

$$\Upsilon = \sqrt{\Pi^{-1}}(I - \Sigma)\tilde{\Upsilon}(I - \Sigma)\sqrt{\Pi^{-1}}.$$

PROOF.  The proof is similar to Theorem 5.5, and therefore only an outline of the proof is given.

Following the proof of (5.11), there exists $\mathbf{s}^\nu$ lying between $\bar{\mathbf{T}}_\infty^\nu$ and $\boldsymbol{\tau}$ such that

$$(\bar{\mathbf{T}}_\infty^\nu - \boldsymbol{\tau})\sqrt{N^\nu\Pi}U(\mathbf{s}^\nu) = \mathbf{Y}^\nu(\bar{\mathbf{T}}_\infty^\nu + \boldsymbol{\zeta}^\nu) + \mathbf{A}_\nu + (\boldsymbol{\zeta}^\nu - \boldsymbol{\zeta})\sqrt{N^\nu\Pi},$$

where $U(\cdot)$ is defined in Theorem 5.5, and

$$\mathbf{A}_\nu = \{\mathbf{r}^\nu(\bar{\mathbf{T}}_\infty^\nu + \boldsymbol{\zeta}^\nu)P^\nu - \mathbf{r}(\bar{\mathbf{T}}_\infty^\nu + \boldsymbol{\zeta}^\nu)\}\sqrt{N^\nu\Pi}.$$

Since $\mathbf{r}(\mathbf{t}) = \mathbf{1} - \mathbf{g}(\mathbf{z}(\mathbf{t}))$, let

$$\begin{aligned}
\tilde{\mathbf{A}}_\nu &= \mathbf{r}(\bar{\mathbf{T}}_\infty^\nu + \boldsymbol{\zeta}^\nu)(\Pi^\nu - \Pi)\sqrt{N^\nu\Pi^{-1}} \\
&= \sqrt{N^\nu}(\boldsymbol{\pi}^\nu - \boldsymbol{\pi})(I - G(\mathbf{z}(\bar{\mathbf{T}}_\infty^\nu + \boldsymbol{\zeta}^\nu)))\sqrt{\Pi^{-1}}.
\end{aligned}$$

It is trivial to show that

$$|\mathbf{A}_\nu - \tilde{\mathbf{A}}_\nu| = |\{\mathbf{r}^\nu(\bar{\mathbf{T}}_\infty^\nu + \boldsymbol{\zeta}^\nu) - \mathbf{r}(\bar{\mathbf{T}}_\infty^\nu + \boldsymbol{\zeta}^\nu)\} P^\nu \sqrt{N^\nu \Pi^{-1}}| \xrightarrow{p} 0 \qquad \text{as } \nu \to \infty.$$

Let $\hat{A}_\nu = \sqrt{N^\nu}\{\boldsymbol{\pi}^\nu - \boldsymbol{\pi}\}(I - \Sigma)\sqrt{\Pi^{-1}}$. Since $G(\mathbf{z}(\bar{\mathbf{T}}_\infty^\nu + \boldsymbol{\zeta}^\nu)) \xrightarrow{p} G(\mathbf{z}(\boldsymbol{\tau} + \boldsymbol{\zeta})) \equiv \Sigma$ [see (5.6)] and $\sqrt{N^\nu}\{\boldsymbol{\pi}^\nu - \boldsymbol{\pi}\} \xrightarrow{D} \mathbf{S}$ as $\nu \to \infty$, it follows that $|\hat{\mathbf{A}}_\nu - \tilde{\mathbf{A}}_\nu| \xrightarrow{p} 0$ as $\nu \to \infty$.

By Lemmas 6.1 and 5.4, and the triangle inequality,

$$|\mathbf{Y}^\nu(\bar{\mathbf{T}}_\infty^\nu + \boldsymbol{\zeta}^\nu) - \tilde{\mathbf{Y}}^\nu(\boldsymbol{\tau} + \boldsymbol{\zeta})| \xrightarrow{p} 0 \qquad \text{as } \nu \to \infty.$$

Therefore

$$(\bar{\mathbf{T}}_\infty^\nu - \boldsymbol{\tau})\sqrt{N^\nu \Pi} = (\tilde{\mathbf{Y}}^\nu(\boldsymbol{\tau} + \boldsymbol{\zeta}) + \hat{\mathbf{A}}_\nu)U(\mathbf{s}^\nu)^{-1} + \mathbf{M}_\nu,$$

where $\mathbf{M}_\nu \xrightarrow{p} \mathbf{0}$ as $\nu \to \infty$.

For all $\nu \geq 1$, $\tilde{\mathbf{Y}}^\nu(\boldsymbol{\tau} + \boldsymbol{\zeta})$ and $\hat{\mathbf{A}}_\nu$ are independent, with $\tilde{\mathbf{Y}}^\nu(\boldsymbol{\tau} + \boldsymbol{\zeta}) \xrightarrow{D} N(\mathbf{0}, \Xi)$ and $\hat{\mathbf{A}}_\nu \xrightarrow{D} \mathbf{S}(I - \Sigma)\sqrt{\Pi^{-1}}$ as $\nu \to \infty$. Therefore Theorem 6.2 follows since $U(\mathbf{s}^\nu) \xrightarrow{p} U$ as $\nu \to \infty$.  $\square$

The adaption of Theorem 6.2 to other (random) group allocation schema is trivial. Theorem 5.5 is the special (deterministic) case where $\mathbf{S} \equiv \mathbf{0}$.

**7. Special cases.** The generic results developed in the preceding sections can be specialized to a number of specific examples. These highlight the wide-ranging applicability of the model described in Section 2.

7.1. *Ball–Clancy* (1993) *model.* This very general model is described in full detail in [6]. The important features from our perspective are the following. There are $m$ groups, and for $\nu \geq 1$, there are $\mathbf{a}^\nu$ initial infectives and $\mathbf{N}^\nu$ initial susceptibles. The susceptibles remain fixed within their initial group. However, once infectious, an individual can move between the $m$ groups. For $1 \leq i \leq m$, let $\mathbf{I}^i = (I_1^i, I_2^i, \ldots, I_m^i)$, where for $1 \leq k \leq m$, $I_k^i$ denotes the time spent in group $k$ by an infective who originates from group $i$, with $I^i = \sum_{k=1}^m I_k^i$. For $1 \leq i, j, k \leq m$, let $\frac{1}{N^\nu}\beta_{kj}^{(i)}$ denote the rate at which an infective, who originates from group $i$, and is currently in group $j$, makes infectious contacts with a given individual in group $k$. Let $B_i = (\beta_{kj}^{(i)})$ $(1 \leq i \leq m)$. Clearly for $1 \leq i, k \leq m$, $V_{i,k}^\nu = 1 - \exp(-\frac{1}{N^\nu}\sum_{j=1}^m \beta_{kj}^{(i)} I_j^i)$, and hence, $N_k^\nu V_{i,k}^\nu \approx \pi_k \sum_{j=1}^m \beta_{kj}^{(i)} I_j^i$. Thus for $1 \leq i \leq m$,

$$\mathbf{V}_i^\nu N^\nu \Pi^\nu \xrightarrow{D} \mathbf{U}_i = \mathbf{I}^i B_i^T \Pi \qquad \text{as } \nu \to \infty.$$

Therefore, if for $1 \leq i \leq m$, $\mathbb{E}[I^i] < \infty$ and $\text{var}(I^i) < \infty$ (cf. [6], Section 4, page 727), it is straightforward to show that the conditions stated in Sections 4 and 5 are satisfied.

7.2. *Ball–Clancy* (1995) *model.* In [7] an epidemic in a homogeneously mixing population with $m$ different types of infectives is considered. An individual on becoming infected chooses its type at random from $\{1, 2, \ldots, m\}$ according to a probability distribution $\boldsymbol{\pi}$, irrespective of the type of the infector. An equivalent construction is to randomly allocate a type to each individual according to $\boldsymbol{\pi}$ prior to the epidemic. Then should an individual become infected, it becomes an infective of its allocated type. Therefore this is a special case of the model described in Section 6.

7.3. *Random graph models.* The main motivation for the current work is the study of the final size distribution of epidemics upon random graphs. A population (and epidemic) upon a random graph $\mathcal{G}$ are defined as follows. Two individuals, $i$ and $j$, are said to be acquaintances, if in the graph $\mathcal{G}$ an edge exists between vertices $i$ and $j$.

A wide variety of extensions of the Bernoulli random graph can be shown to satisfy the randomized Reed–Frost criterion. We shall illustrate this range by considering static and dynamic multitype Bernoulli random graphs. The key property that we shall require is that the existence of an edge between any two vertices is independent of the remainder of the graph. This will allow the epidemic and random graph to be constructed in unison; see [16].

For a static multitype Bernoulli random graph $\mathcal{G}$, for $1 \leq i, j \leq m$, let $\frac{\alpha_{ij}}{N^\nu}$ $(\alpha_{ij} \geq 0; N^\nu \geq \max\{\alpha_{ij}\})$ denote the probability that there exists an edge between a given vertex of type $i$ and a given vertex of type $j$. We shall require that the random graph is irreducible, that is, there is a nonzero probability of there existing a pathway of edges between any pair of vertices. Let $W_{i,j}$ denote the probability that a type-$i$ individual, while infectious, makes an infectious contact with a type-$j$ acquaintance. Thus $V_{i,j}^\nu = \frac{\alpha_{ij}}{N^\nu} W_{i,j}$ is the probability that a type-$i$ infective, while infectious, makes an infectious contact with a type-$j$ susceptible. Therefore the conditions of Sections 4 and 5 are satisfied. The underlying multitype random graph can then be constructed in a similar manner to [16] with an edge existing between two vertices $k$ and $l$, if in the epidemic individual $k$ infects $l$ or vice versa. An edge exists between the vertices $k$ and $l$, according to the correct conditional distribution, if neither individual $k$ nor $l$ infects the other individual.

An interesting special case of the static multitype Bernoulli random graph with random type allocation is the mixed Bernoulli random graph. The mixed Bernoulli random graph epidemic model assumes that the infectivity, and connectivity, of all individuals are independent and identically distributed according to random variables $W$ and $D$, respectively, with $W$ and $D$ independent. That is, for individuals $k$ and $l$ with connectivity random variables $D_k$ and $D_l$, respectively, the probability that $k$ and $l$ are acquaintances is proportional to $D_k D_l$. Then given they are acquaintances, the probability that individual $k$, while infectious, makes infectious contacts with individual $l$ is $W_k$. We assume that there exist $m \in \mathbb{N}$, and $\boldsymbol{\theta}, \boldsymbol{\pi} \in \mathbb{R}^m$, such that $\sum_{i=1}^{m} \pi_i = 1$ and $\mathbb{P}(D = \theta_j) = \pi_j$ $(1 \le j \le m)$. Therefore, independently at random assign individuals to one of $m$ groups according to $\boldsymbol{\pi}$, that is, according to $D_i$. Then for $1 \le i, j \le m$, let an individual of type $i$ and an individual of type $j$ be acquaintances with probability $\theta_i \theta_j / N^\nu$. Thus $V_{i,j}^\nu \overset{D}{=} \theta_i \theta_j W / N^\nu$ and the conditions of Section 6 are satisfied.

Epidemics upon random graphs with prescribed degree have received considerable attention in the literature; see, for example, [3, 4, 17]. That is where a vertex, $i$ say, has $E_i$ edges connecting the vertex to $E_i$ distinct vertices. The random variables $\{E_1, E_2, \ldots, E_{N^\nu}\}$ are assumed to be (asymptotically) independent and identically distributed according to a nonnegative, integer-valued random variable $E$. Particular interest has been devoted to the so-called scale-free graphs [17], where the mean of $E$ exists, but the variance does not. Unfortunately, such epidemic models do not permit a (multitype) randomized Reed–Frost representation since the existence of an edge between any pair of vertices is no longer independent of the remainder of the graph. However, in most real-life scenarios it is realistic to assume that $\mathbb{E}[E] \le \mathrm{var}(E)$, and so, the mixed Bernoulli model with $D$ chosen such that $\mathbb{E}[E] = \mathbb{E}[D]^2$ and $\mathrm{var}(E) = \mathbb{E}[D]^2(1 + \mathrm{var}(D))$ will often provide a suitable alternative. (We have a scale-free mixed Bernoulli random graph if the mean of $D$ exists, but the variance does not.) Thus it would be interesting to extend the results for the mixed Bernoulli random graph to continuous and/or scale-free $D$.

Finally, we consider dynamic Bernoulli random graph epidemics, which have previously been considered in [1] and [4]. For a random graph $\mathcal{G}$, and for $1 \le i, j \le m$, consider two distinct vertices, one of type $i$, vertex $k$, say, and one of type $j$, vertex $l$, say. In the dynamic random graph, we alternate between there existing an edge between the vertices $k$ and $l$, and there being no edge between the vertices $k$ and $l$. The time from the formation of an edge between vertices $k$ and $l$, until the dissolution of the edge is distributed according to $S_{ij}^- \sim \mathrm{Exp}(\rho_{ij}^-)$. The time from the dissolution of the edge between vertices $k$ and $l$ until the reforming of the edge is distributed according to $N^\nu S_{i,j}^+$, where $S_{i,j}^+ \sim \mathrm{Exp}(\rho_{ij}^+)$. We shall assume that for $1 \le i, j \le m$, $\rho_{ij}^+, \rho_{ij}^- > 0$, that is, edges can form and dissolve between all types of vertices. Thus for $t \ge 0$, the two individuals corresponding to vertices $k$ and $l$ are acquaintances at time $t$, if an edge exists between vertices $k$ and $l$ at time $t$. Therefore, in equilibrium $\alpha_{ij}^\nu = \frac{\rho_{ij}^+}{\rho_{ij}^+ + N^\nu \rho_{ij}^-}$ is the probability that the two individuals are acquainted at any given point in time. We shall assume that the dynamic Bernoulli random graph is in equilibrium.

Only minor modifications are required for the following cases: no edges exist between certain types of vertices (all that we require is the graph $\mathcal{G}$ is again irreducible), the graph is part static and part dynamic and general probability distributions for $S_{i,j}^-$ and $S_{i,j}^+$.

For simplicity in exposition we shall consider a generalization of the model considered in [4]. For $1 \le i \le m$, infected individuals of type $i$ have independent and identically distributed lifetimes according to $Q_i$, where $\mathbb{E}[Q_i^2] < \infty$. While infectious (and acquainted), an individual of type $i$ makes infectious contacts with an acquaintance of type $j$ at the points of a homogeneous Poisson point process with rate $\beta_{ij}$. As noted in [4], the probability of more than one partnership between two individuals during the course of an infectious period is negligible. We therefore have the following two cases to consider: first, where the two individuals are acquaintances at the beginning of the infectious period, which we shall term initial acquaintanceships, and second, where the two individuals become acquainted during the infectious period, which we shall term secondary acquaintanceships. Suppose that an individual of type $i$, and an individual of type $j$, are acquaintances at the point of time, $t$ say, at which individual $i$ becomes infected. Suppose that the individuals remain acquainted until time $t + Z_{i,j}$, where $Z_{i,j} \overset{D}{=} S_{i,j}^-$. Conditioning upon $Q_i = q_i$, the probability that an initial acquaintance of type $j$ is infected is $1 - \exp(-\beta_{ij}(q_i \wedge Z_{i,j}))$. The probability that during the infective's infectious period a link forms between the infective and an initially unacquainted individual of type $j$ is $1 - \exp(-\frac{1}{N^\nu}\rho_{ij}^+ q_i) \approx \frac{1}{N^\nu}\rho_{ij}^+ q_i$. Conditional upon a link between the two individuals being formed, the time (relative to the start of the infectious period) at which the acquaintanceship is formed has probability density function $g(x) = \frac{\rho_{ij}^+}{N^\nu}\exp(-x\rho_{ij}^+/N^\nu)/\{1 - \exp(-q_i\rho_{ij}^+/N^\nu)\} \approx 1/q_i$ ($0 \le x \le q_i$), and $g(x) = 0$ (otherwise). Therefore the probability of infecting an individual of type $j$ who is not an initial acquaintance is approximately

$$\frac{1}{N^\nu}\rho_{ij}^+ q_i \int_0^{q_i} 1/q_i\{1 - \exp\{-\beta_{ij}((q_i - x) \wedge S_{ij}^-)\}\}\,dx.$$

Hence, $V_{i,j}^\nu = M_{i,j}^\nu + o(1/N^\nu)$, where

$$\begin{aligned}
M_{i,j}^\nu &= \alpha_{ij}^\nu\{1 - \exp(-\beta_{ij}(Q_i \wedge Z_{i,j}))\} \\
&\quad + \frac{\rho_{ij}^+}{N^\nu}\int_0^{Q_i}\{1 - \exp(-\beta_{ij}((Q_i - x) \wedge S_{i,j}^-))\} \\
&= \alpha_{ij}^\nu\frac{\beta_{ij}}{\rho_{ij}^- + \beta_{ij}}\{1 - \exp(-(\beta_{ij} + \rho_{ij}^-)Q_i)\} \\
&\quad + \frac{\rho_{ij}^+}{N^\nu}\frac{\beta_{ij}}{\rho_{ij}^- + \beta_{ij}}\Big\{Q_i - \frac{1}{\rho_{ij}^- + \beta_{ij}}(1 - \exp(-Q_i(\rho_{ij}^- + \beta_{ij})))\Big\}.
\end{aligned}$$

Therefore, for $1 \le i \le m$, there exists a random vector $\mathbf{U}_i$ with finite mean and covariance matrix such that $\mathbf{V}_i^\nu N^\nu \Pi^\nu \overset{D}{\longrightarrow} \mathbf{U}_i$.

**Acknowledgments.** I would like to thank the referee and an Associate Editor for their helpful comments. Their comments have improved the presentation of the paper and provided the inspiration for Section 6.

School of Mathematics
University of Manchester
Sackville Street
Manchester
M60 1QD
United Kingdom
E-mail: P.Neal-2@manchester.ac.uk